\documentclass[a4paper,fleqn]{cas-dc}

\usepackage[authoryear,longnamesfirst]{natbib}
\usepackage{amsmath,amssymb,bm}
\usepackage{graphicx}
\usepackage{xcolor}
\usepackage{booktabs}
\usepackage{subcaption} % Existing subfigure environments in the appendix.
\usepackage{tikz}       % Existing TikZ diagram in the FEM appendix.
\usepackage{ragged2e}   % Existing \justifying figure captions.
\usepackage{hyperref}   % Existing hyperlinks to run diagnostics.
\usepackage[section]{placeins} % Keep floats within their sections where possible.
\usepackage{float} % Enable [H] placement for single-column figures.

\newcommand{\R}{\mathbb{R}}
\newcommand{\wphi}{\omega_{\phi}}
\newcommand{\wchi}{\omega_{\chi}}

\newcommand{\RMSx}[1]{\mathrm{RMS}_x\!\left[#1\right]}
\DeclareMathOperator{\sech}{sech}

\begin{document}
\let\WriteBookmarks\relax

\shorttitle{Metastability Landscapes in Nonlinear Ghost-Normal Scalar Fields}
\shortauthors{Wysong et al.}
\title[mode=title]{Metastability Landscapes in Nonlinear Ghost-Normal Scalar Fields}

% Authors' ORCID IDs and e-mail addresses transferred from REVTeX metadata.
% No corresponding author is designated in the source. Add \cormark/\cortext
% only after the authors confirm who should be the corresponding author.
\author[1]{Jax Wysong}[orcid=0000-0001-5766-4332]
\ead{wysongj2@msu.edu}

\author[2]{Samara Overvaag}[orcid=0000-0001-7123-7003]
\ead{samara.overvaag@jacks.sdstate.edu}

\author[3]{Hyun Lim}[orcid=0000-0002-8435-9533]
\ead{hyunlim@lanl.gov}

\author[2]{Jung-Han Kimn}
\ead{jung-han.kimn@sdstate.edu}

\affiliation[1]{organization={Department of Computational Mathematics, Science, and Engineering, Michigan State University},
    city={East Lansing},
    postcode={48823},
    state={MI},
    country={USA}}
    
\affiliation[2]{organization={Department of Mathematics and Statistics, South Dakota State University},
    city={Brookings},
    postcode={57007},
    state={SD},
    country={USA}}

\affiliation[3]{organization={Applied Computer Science (CAI-1) and Center for Theoretical Astrophysics, Los Alamos National Laboratory},
    city={Los Alamos},
    postcode={87545},
    state={NM},
    country={USA}}

% No author contributions, title footnotes, or funding notes supplied in source.
% Insert \credit{...}, \tnotemark/\tnotetext, or \nonumnote if needed.

\begin{abstract}
We investigate finite time dynamics in nonlinear classical scalar field systems containing one ghost and one normal kinetic sector. Using a spacetime finite element formulation, we study the dependence of threshold crossing and numerical continuation times on initial amplitude, characteristic wavenumber, spectral tilt, phase structure, and spatial localization. In the quartic model, larger amplitudes generally lead to earlier growth, while higher wavenumbers and ultraviolet-biased spectra favor longer resolved bounded evolution. Smooth compact initial data exhibit a broad region with no observed growth threshold crossing over the simulated interval. 
A direct comparison of two Gaussian evolutions shows that enhanced low-wavenumber power can precede rapid field growth, although this observation does not establish a universal spectral-transfer mechanism. We also compare the quartic interaction with a lifted $\phi^6$ potential. The latter gives a markedly different continuation landscape, with a broad long-converged region followed by shorter accessible terminal times at larger amplitudes. Localization diagnostics do not identify a persistent oscillon core in the selected cases. By distinguishing amplitude growth, spectral reorganization, solver-limited termination, and finite-time censoring, our results characterize a nontrivial metastability landscape without implying global or quantum stability.
\end{abstract}

% Optional Elsevier front matter; drafted from claims in the original manuscript.

\maketitle

\section{Introduction}
\label{sec:intro}
The absence of ghost-like excitations is traditionally regarded as a fundamental prerequisite for the internal consistency of any physical theory. The presence of negative kinetic terms, corresponding to dynamical modes carrying negative norm states or yielding Hamiltonians unbounded from below, typically signals severe instabilities at both the classical and quantum levels. At the quantum level, ghosts are associated with violations of unitarity or catastrophic vacuum decay through the uncontrolled production of positive–negative energy pairs~\cite{Cline:2003gs,Woodard:2006nt}. At the classical level, they lead to runaway solutions due to the absence of a lower bound on the energy functional. For these reasons, ghost degrees of freedom are usually regarded as pathological and are excluded in standard model building.

Nevertheless, ghost-like fields arise rather naturally in a wide range of physical frameworks. In cosmology, effective descriptions of dark energy with equation of state parameter $w<-1$, so called phantom models, are most simply realized by scalar fields with negative kinetic terms~\cite{Caldwell:2002ii,Cline:2003gs}. In modified theories of gravity, higher-curvature corrections introduce additional propagating degrees of freedom, including massive spin-2 ghosts in quadratic gravity~\cite{Stelle:1977ry,Stelle:1978yy, salvio19}. Similarly, higher derivative scalar field theories and effective field theories (EFTs) containing operators such as $(\Box \phi)^2$ or $R^2$ often propagate extra modes whose stability depends delicately on the structure of the theory~\cite{Woodard:2006nt,Simon:1990ic}. The appearance of such modes reflects a deeper tension between theoretical consistency, bounded Hamiltonians, unitarity, well-posedness, and the desire to incorporate ultraviolet corrections or explain observational phenomena within a unified framework.

A well-known result governing higher derivative systems is the Ostrogradsky theorem~\cite{Ostrogradsky:1850fid}, which states that non-degenerate Lagrangians with higher time derivatives lead to Hamiltonians that are linear in at least one canonical momentum and therefore unbounded from below. This mathematical result is often interpreted as implying unavoidable ghost instabilities in generic higher derivative theories. However, the scope of the theorem is precise: it applies to non-degenerate systems. Degenerate higher derivative theories such as Galileons and more general degenerate higher order scalar-tensor theories can evade Ostrogradsky instabilities by imposing primary constraints that eliminate the would-be ghost degree of freedom~\cite{Nicolis:2008in,Langlois:2015cwa}. These constructions demonstrate that higher derivatives do not automatically imply propagating ghosts, and they play a central role in modern scalar–tensor model building.

From the perspective of EFTs, higher derivative operators arise naturally as suppressed corrections in an expansion in powers of $E/\Lambda$, where $\Lambda$ is the cutoff scale~\cite{Burgess:2007pt}. In this context, ghost-like excitations often appear at or above the cutoff and are interpreted as artifacts of truncating the EFT rather than genuine propagating degrees of freedom within its regime of validity~\cite{Simon:1990ic,Donoghue:1994dn}. When treated perturbatively, so that higher-derivative terms are viewed as small corrections, runaway solutions associated with Ostrogradsky modes can lie outside the domain of validity of the EFT. This viewpoint suggests that the catastrophic instabilities predicted by a strict application of the theorem may not always manifest physically, particularly when the theory is interpreted as an effective description valid only below some scale.

Recent work has explored several interesting directions that further nuance the conventional wisdom about ghosts. For instance, certain higher-derivative scalar theories allow for a stable crossing of the phantom divide ($w=-1$) using a single degree of freedom without generating ghost-like perturbations around cosmological backgrounds~\cite{Vikman:2004dc,Creminelli:2006xe}. While general higher-derivative theories are subject to Ostrogradsky instabilities, carefully structured interactions or degeneracy conditions can eliminate the dangerous mode at the perturbative level.

Moreover, although instabilities associated with Ostrogradsky ghosts are derived rigorously at the mathematical level, recent investigations indicate that catastrophic behavior need not be generic in physically relevant settings. In the context of classical point-particle models, the existence of additional integrals of motion has enabled proofs of global stability within a class of higher-derivative systems~\cite{Deffayet2022,Deffayet2023}. In these cases, all phase-space trajectories remain bounded for arbitrary initial conditions, providing explicit counterexamples to the expectation that Ostrogradsky systems must exhibit runaway behavior. See also~\cite{Diez2025} for a systematic discussion and classification of stability results in higher-derivative mechanical systems.

Another arena where progress has been made is classical field theory, particularly in effective field theories of gravity that explicitly include higher-derivative operators. Quadratic gravity, for example, propagates a massive spin-2 ghost at the linearized level~\cite{Stelle:1977ry}, yet recent work has investigated the well-posedness of its initial value formulation and the possibility of controlled evolutions~\cite{Noakes:1983,Held2021,Held2023,Held2025}. These studies suggest that, under suitable conditions and within certain regimes, the dynamics can remain well-defined for long times, at least at the classical level. Furthermore,~\cite{Deffayet2025} presents long-lived time evolutions in a system of two scalar fields, one normal and one ghost, coupled through non-derivative interactions in $(1+1)$-dimensional Minkowski spacetime. 

Importantly, recent analytical results on nonlinear wave and Klein–Gordon systems have established small-data global stability in a class of models that do not rely on positive-definite energy functionals~\cite{Held2025b}. In these systems, sufficiently small initial data lead to globally regular solutions whose dispersive decay suppresses nonlinear growth, even in the presence of ghost-like interactions. This perspective suggests that instability is not an automatic consequence of the Hamiltonian structure, but rather a dynamical phenomenon controlled by the interplay between nonlinearity, dispersion, and initial data.
These developments motivate a more refined view: classical ghost dynamics should be organized not by a binary stable/unstable classification, but by a hierarchy of stability notions.

A remaining problem is to determine how much of this long-lived behavior is organized by the structure of the initial data rather than by one special waveform or one particular definition of numerical failure.
The nonlinear field equations admit families with very different spatial and spectral characteristics: coherent plane waves, localized packets, broadband random spectra, and initially compact perturbations. Such families need not produce the same finite-time behavior even when a nominal field-amplitude scale is held fixed. Moreover, when the equations are solved by an iterative nonlinear method, loss of algebraic convergence is not automatically a physical instability. These considerations motivate both the parameter-space surveys and the explicit separation of physical diagnostics from numerical continuation
times used below.

We examine two coupled scalar fields with opposite-sign kinetic sectors and use the quartic interaction $V_{22}=\lambda_{22}\phi^2\chi^2$ as our baseline model. We first map operational growth-threshold times for plane waves and Gaussian packets, then vary the initial spectral tilt at fixed field root-mean-square amplitude in random-phase ensembles. Additional phase-correlated and smooth compact-support data probe the roles of relative waveform alignment and initial localization. Compact support here refers to the initial profile on a periodic computational domain, and therefore should not be confused with scattering on unbounded space.

The spectral surveys motivate a second question: is long-wavelength growth merely correlated with the initial spectrum, or can it emerge dynamically before rapid amplitude growth? We address this by comparing time-resolved Fourier-power distributions and field amplitudes in two Gaussian runs.
Fourier power is an amplitude diagnostic, not the signed energy of the normal and ghost sectors, and we do not interpret one such comparison as proof of a universal inverse energy cascade. We additionally consider a lifted $\phi^6$ interaction motivated by models that support localized oscillatory configurations in positive-energy scalar theory~\cite{GleiserSicilia2009,Amin2010,HindmarshSalmi2008}.
In the present implementation this interaction also modifies the linearized mass spectrum. The resulting comparison therefore probes sensitivity to the complete nonlinear field-theory model rather than isolating a sextic stabilization mechanism.

This work can be described as threefold. First, we establish systematic finite-time maps across several initial data families, including the spectral-tilt ordering and its phase-realization variability. Second, we identify a representative increase in low-wavenumber power preceding rapid field growth. Third, we distinguish measured operational threshold crossings, right-censored evolutions, and solver-limited terminal times, which is
particularly important for interpreting the higher order potential runs. These observations characterize a nontrivial metastability landscape but do not constitute a theorem of global stability or a demonstration of quantum consistency.

The paper is organized as follows.  Sec.~\ref{sec:model} introduces the models, potentials, initial data, numerical methods, and metastability diagnostics. Sec.~\ref{sec:quartic_landscapes} maps the lifetime landscape of the quartic ghost-normal model. Sec.~\ref{sec:spectral_mechanism} analyzes spectral-transfer and phase-coherence mechanisms.  Sec.~\ref{sec:phi6} studies the lifted $\phi^6$ model and oscillon-mediated metastability.  Sec.~\ref{sec:conclusion} summarizes implications, limitations, and open problems.  Throughout this work, we employ natural units with $c=1$.

\section{Models, observables, and numerical setup}
\label{sec:model}
We consider a system with a Lagrangian,
\begin{align}
\mathcal{L} ={}& \frac{1}{2}\phi(\Box+m_{\phi}^2)\phi \nonumber\\
 &+ \frac{\gamma}{2}\chi(\Box+m_{\chi}^2)\chi
   + V(\phi,\chi).
\label{eq: Ghost Lagrangian}
\end{align}
Where $\Box = \partial_t^2 - \nabla^2$ is the D'Alembertian operator, $\nabla$ is a gradient, $\phi$ and $\chi$ are scalar fields, and $V(\phi, \chi)$ is a nonlinear potential that couples the scalar fields. The ghost appears in this Lagrangian if the scalar fields are assigned oppositely signed kinetic energy terms, $\gamma = -1$.

We find the strong equations of motion
\begin{align}
\Box \phi + m_{\phi}^2 \phi + \partial_\phi V &= 0, \label{eq:phi_strong_EOM} \\
\Box \chi + m_{\chi}^2 \chi + \gamma \partial_\chi V &= 0.\label{eq:chi_strong_EOM}
\end{align}

As described in~\cite{Deffayet2025}, the potential can be split according to
\begin{align}
    V_\phi(\phi) &\equiv V(\phi,\chi=0), \\ 
    V_\chi(\chi) &\equiv V(\phi=0,\chi), \\ 
    V_\mathrm{int}(\phi, \chi) &= V(\phi,\chi) - V_\phi (\phi) - V_\chi (\chi), 
    \label{eqn:potentials}
\end{align}
where $V_\phi$ and $V_\chi$ are considered to be self-interaction potentials while $V_\mathrm{int}$ is the ghostly interaction potential. 
We define energies by integrating Hamiltonian density over the spatial domain
\begin{align}
    H &= \int_\Omega \mathcal{H}, \\
    H_\phi &= \int_\Omega \mathcal{H}_\phi, \\
    H_\chi &= \int_\Omega \mathcal{H}_\chi, \\
    H_\mathrm{int} &= \int_\Omega \mathcal{H}_\mathrm{int}, 
\end{align}
with the following Hamiltonian densities
\begin{align}
\mathcal{H} &= \frac12\big[(\partial_t\phi)^2+(\nabla\phi)^2
                 +m_\phi^2\phi^2\big] \nonumber\\
 &\quad+\frac\gamma2\big[(\partial_t\chi)^2+(\nabla\chi)^2
                 +m_\chi^2\chi^2\big] \nonumber\\
 &\quad+V(\phi,\chi),\\
\mathcal{H}_\phi &= \frac12\big[(\partial_t\phi)^2+(\nabla\phi)^2
                     \nonumber\\
 &\qquad+m_\phi^2\phi^2\big]+V_\phi(\phi),\\
\mathcal{H}_\chi &= \frac\gamma2\big[(\partial_t\chi)^2+(\nabla\chi)^2
                     \nonumber\\
 &\qquad+m_\chi^2\chi^2\big]+V_\chi(\chi),\\
\mathcal{H}_\mathrm{int} &= \mathcal{H}-\mathcal{H}_\phi-\mathcal{H}_\chi.
\end{align}
Since $H$ corresponds to the total energy, it should be conserved regardless of the sign of $\gamma$. 

Numeric evidence in Ref.~\cite{Deffayet2025} indicates that the quartic interaction can yield delayed or benign classical ghost behavior rather than instantaneous runaway.  We use this model as the baseline for the metastability maps below. Representative diagnostic examples are discussed in
Appendix~\ref{sec:appx:runaway-diagnostics}.

\subsection{Potentials and organizing questions}
\label{ssec:potentials}

Unless otherwise stated, we use the polynomial interaction
\begin{equation}
    V_{22}=\lambda_{22}\phi^2\chi^2 .
\end{equation}
This model is used as a baseline for mapping how metastability depends on amplitude, characteristic wavenumber, spectral tilt, localization, and phase coherence. 
For this potential,
\begin{equation}
V_\phi=V_\chi=0,
\qquad
V_{\rm int}=V_{22},
\end{equation}
and
\begin{equation}
\partial_\phi V_{22}
=2\lambda_{22}\phi\chi^2,
\qquad
\partial_\chi V_{22}
=2\lambda_{22}\phi^2\chi.
\end{equation}

We additionally study a lifted polynomial interaction
of the form
\begin{equation}
U(W)=\frac{1}{2} m_6^2W
-\frac{\lambda}{4}W^2
+\frac{g}{6}W^3,
\qquad
W=\phi^2+\chi^2.
\label{eq:lifted_potential}
\end{equation}
The lifted model is analyzed separately in
Sec.~\ref{sec:phi6}. In the implemented equations,
its quadratic contribution is added to the bare
mass terms; consequently, the two models also
differ in their linearized mass spectra.

Our comparison therefore concerns the finite-time
dynamics of two specified interaction models, rather
than isolating the effect of the sextic term at
fixed linearized masses.

\subsection{Numerical Methods}
\label{sec:methods}

We discretize the spatial and temporal domains using the finite element method (FEM) on spacetime. 
This method does not require time stepping. Instead, the entire spacetime slab is solved simultaneously. The use of a spacetime FEM offers several advantages over more conventional finite-difference schemes with explicit time integration, particularly for nonlinear systems that may develop steep gradients or near-singular features. In the spacetime formulation, spatial and temporal discretizations are treated on equal footing, allowing the solution to be obtained over an entire spacetime slab without the need for iterative time stepping. This approach improves global consistency and avoids cumulative integration errors that can accumulate in explicit schemes. Moreover, the variational structure of the FEM preserves key conservation properties such as total energy more accurately under coarse resolution. 

The method also accommodates flexible boundary conditions and naturally incorporates adaptive refinement strategies in both space and time, which is essential for capturing localized structures such as sharp ghost-induced instabilities. Compared with explicit finite-difference methods, the space–time FEM provides enhanced numerical stability for stiff or mixed-sign kinetic systems, where traditional Courant–Friedrichs–Lewy (CFL) limits might become restrictive. 

These properties make the method particularly suitable for exploring nonlinear dynamics in ghost-containing field theories, where stability and long-time accuracy are critical. The Portable, Extensible Toolkit for Scientific Computation (PETSc)~\cite{petsc-web-page} library  is used to complement our base C code. PETSc has a host of convenient options for parallel vector/matrix operations, as well as a suite of different linear (KSP) and nonlinear (SNES) solvers. 
A more detailed discussion of spacetime FEM including element stiffness matrix calculations is provided in Appendix~\ref{sec:appx:sFEM}.

The nonlinear algebraic problem is solved with PETSc SNES. For the longer calculations we use a warm-start continuation procedure: after converging a solution on a time domain $[0,T]$, the converged solution is used to construct the initial Newton iterate for an extended domain $[0,T+\Delta T]$.
The procedure can increase the attainable terminal time by improving the initial nonlinear guess. The last terminal time reached by this sequence, $T_{\rm conv}$, is therefore a diagnostic of numerical continuation, not an independently measured physical instability time. In particular, a change in $T_{\rm conv}$ upon changing the initial Newton iterate must not be interpreted as a physical stabilization effect.

\subsection{Initial Data}
\label{sec:id}
To investigate how finite-time dynamics depends on the initial configuration, we consider six families of initial data: plane waves, Gaussian packets, random-phase broadband spectra, phase-correlated plane waves, localized $\sech$ profiles, and smooth compactly supported perturbations.

\subsubsection{Plane waves}
\label{ssec:id_plane}
The plane wave initial conditions are given by
\begin{align}
    \phi_0(x) &= A_\phi \, \sin \bigl( k_\phi (x - x_\phi) \bigr), \nonumber \\
    \dot{\phi}_0(x) &= - c_\phi \, k_\phi \, A_\phi \, \cos \bigl( k_\phi (x - x_\phi) \bigr), \nonumber\\
    \chi_0(x) &= A_\chi \, \sin \bigl( k_\chi (x - x_\chi) \bigr), \nonumber\\
    \dot{\chi}_0(x) &= - c_\chi \, k_\chi \, A_\chi \, \cos \bigl( k_\chi (x - x_\chi) \bigr). \nonumber
\end{align}
where
\begin{equation}
\begin{aligned}
c_\phi &= \pm \sqrt{\frac{k_\phi^2 + m_\phi^2}{k_\phi^2}}, \\
c_\chi &= \pm \sqrt{\frac{k_\chi^2 + m_\chi^2}{k_\chi^2}}.
\end{aligned}
\end{equation}
We fix $x_\phi = 0$, $x_\chi = L/3$, $A_\phi = A_\chi = A$, $k_\phi = k_\chi/2 = k$, and choose the two plane waves to move in opposite directions.

$L$ is the spatial length that we implement in the simulation and is always set to 1.0.
Despite maintaining $L = 1.0$ we explicitly include it in many formulae to keep units less ambiguous. 

\subsubsection{Gaussian packets}
The Gaussian wave packet initial data are defined by
\begin{align*}
\phi_0(x) &= A_{\phi}\exp\!\left[-\frac{(x-x_{\phi})^2}{2\ell_{\phi}^2}\right],\\
\dot\phi_0(x) &= \frac{A_{\phi}c_{\phi}(x-x_{\phi})}{\ell_{\phi}^2}
                 \\[-2pt]
 &\qquad\times\exp\!\left[-\frac{(x-x_{\phi})^2}{2\ell_{\phi}^2}\right],\\
\chi_0(x) &= A_{\chi}\exp\!\left[-\frac{(x-x_{\chi})^2}{2\ell_{\chi}^2}\right],\\
\dot\chi_0(x) &= \frac{A_{\chi}c_{\chi}(x-x_{\chi})}{\ell_{\chi}^2}
                  \\[-2pt]
 &\qquad\times\exp\!\left[-\frac{(x-x_{\chi})^2}{2\ell_{\chi}^2}\right].
\end{align*}
where $c_{\phi}$ and $c_{\chi}$ are defined oppositely (1 and -1) so that the waves travel in opposite directions or defined equivalently (1 and 1) for co-moving propagation. The initial positions are described by $x_{\phi} = 0.3L$ and $x_{\chi} = 0.7L$ and the amplitudes and lengths are kept consistent throughout, $A \equiv A_{\phi} = A_{\chi}$ and $\ell \equiv \ell_{\phi} = \ell_{\chi}$. Also, we will use the relations $k = 1/(4\ell)$ and $k \times L/(2\pi) = C$, where $C$ will be the quantity that is varied in the simulations.

\subsubsection{Colored-noise spectra (IR/UV-tilted)}
\label{ssec:id_colored}
Choose integers $n_1\le n\le n_2$, phases $\theta_n\in[0,2\pi)$, and a spectral tilt $n_s\in\R$. Define
\begin{align}
  \phi(x,0) &= \mathcal C_\phi \sum_{n=n_1}^{n_2} k_n^{\,n_s/2}\,\cos\!\big(k_n x+\theta_n\big),\\
  \dot\phi(x,0) &= -\sum_{n=n_1}^{n_2} s_n\frac{\wphi(k_n)}{k_n}
                       \nonumber\\
 &\quad\times\partial_x\!\Big[\mathcal C_\phi k_n^{n_s/2}
          \cos\!\big(k_nx+\theta_n\big)\Big],
\end{align}
where $s_n\in\{+1,-1\}$ selects co-/counter-propagating content. Set
\begin{equation}
  \mathcal C_\phi \;=\; \frac{A}{\sqrt{\tfrac12\sum_{n=n_1}^{n_2} k_n^{\,n_s}}}\,,
\end{equation}
so that $\RMSx{\phi}=A$ (for random phases). Define $(\chi,\dot\chi)$ similarly, optionally with a fixed phase offset $\Delta\phi$.

A tilt $n_s<0$ concentrates power in the infrared (IR), while $n_s>0$ biases toward ultraviolet (UV) modes. This family lets us dial the IR/UV balance at fixed RMS to test a central empirical finding of~\cite{Deffayet2025}: at fixed amplitude, higher-$k$ data are more stable than lower-$k$ data. Beyond the present model, colored-noise starts are canonical in weak wave turbulence and in “integrable turbulence” for the 1D nonlinear Schrödinger equation, where random-phase ensembles seed nonlinear cascades and heavy-tailed statistics~\cite{Nazarenko2011,WalczakPRL2015,RandouxSuret2016}.

\subsubsection{Phase-correlated two-field plane waves}
\label{ssec:id_phase}
For a single carrier $k$,
\begin{align}
  \phi(x,0) &= A\,\cos\!\big(k(x-x_\phi)\big), \\
  \dot\phi(x,0) &= -\frac{\wphi(k)}{k}\,\partial_x \phi(x,0),\\
  \chi(x,0) &= A\,r\,\cos\!\big(k(x-x_\chi)+\Delta\phi\big),\qquad \\
  \dot\chi(x,0) &= -\sigma\,\frac{\wchi(k)}{k}\,\partial_x \chi(x,0),
\end{align}
with amplitude ratio $r>0$, relative phase $\Delta\phi\in[0,\pi]$, and $\sigma\in\{+1,-1\}$ for co-/counter-propagation.

This is the cleanest way to control initial cross-correlation between sectors. Relative phase selects a normal-mode mixture and therefore the direction and rate of early energy transfer, analogous to coupled-oscillator physics and to phase-sensitive transfer in preheating/parametric resonance models \cite{KLS1997,GreeneEtAl1997}.

\subsubsection{Oscillon-like, time-symmetric seeds}
\begin{align}
  \phi(x,0) &= A\,\sech\!\Big(\frac{x-x_0}{\sigma}\Big), & \dot\phi(x,0) &= 0,\\
  \chi(x,0) &= A\,r\,\sech\!\Big(\frac{x-x_0}{\sigma}\Big)\cos\Delta\phi, & \dot\chi(x,0) &= 0.
  \label{IC: oscillon-like}
\end{align}
To add a carrier, we may multiply each profile by $\cos(k_0(x-x_0))$.

Oscillons are long-lived, localized, nearly periodic lumps that appear in many scalar theories \cite{GleiserSicilia2009,Amin2010,HindmarshSalmi2008}. Seeding an oscillon-like profile probes whether the ghostly coupling quenches, destabilizes, or stabilizes such coherent structures, an especially discriminating test of nonlinear self-interaction vs. ghost-driven energy exchange. Lifetimes, frequency drifts, and radiation tails provide sharp quantitative metrics.

\subsubsection{Smooth compact support bumps}
\label{IC:compact-bumps}
The basic profile of smooth compactly supported is chosen as the standard $C^\infty$ bump
function
\begin{equation}
f_R(x)=
\begin{cases}
A\exp\left[-\dfrac{1}{1-(r/R)^2}\right],& |r|<R,\\[6pt]
0,& |r|\ge R.
\end{cases}
\label{eq:compact-bump}
\end{equation}
This function is smooth to all orders and has exactly compact support.  It is
therefore qualitatively different from a Gaussian packet, which is localized but
has nonzero tails everywhere.

We will use the compact bump in three complementary ways.  The first is a
single-packet dispersal test, in which $\phi$ and $\chi$ are initialized with
the same compact support but different amplitudes, phases, or widths.  This
tests whether compact support and dispersion delay the nonlinear transfer into
unstable ghost-normal channels.  The second is a two-packet scattering test,
where initially separated compact packets are given opposite momenta and allowed
to collide.  This probes whether two individually long-lived localized
perturbations remain benign after a nonlinear interaction.  The third is a
matched-data comparison against plane waves and Gaussian packets at fixed
amplitude, spectral centroid, and total initial norm.  This comparison separates
the role of spatial support from the role of spectral support.

We take various forms of initial condition. 
For example, we consider
\begin{align}
\phi(0,x) &= A_\phi f_R(x-x_0), \\
\dot{\phi}(0,x) &= 0,\\
\chi(0,x) &= A_\chi f_R(x-x_0)\cos\Delta\theta, \\
\dot{\chi}(0,x) &= 0,
\end{align}
This tests whether a compact perturbation disperses or transfers power into ghost instability. 
%The parameters $(A,R,r,\Delta\theta,v)$ control the
%amplitude, compact-support radius, ghost-to-normal amplitude ratio, relative phase, and packet velocity.  
%The key observables for these runs are the lifetime $T_{\rm life}$, the localization radius $\sigma_x(t)$ or $\sigma_r(t)$, the infrared-to-ultraviolet spectral ratio $R_{\rm IR/UV}(t)$, and the ghost-normal coherence
$\mathcal{C}_{\phi\chi}(t)$.  
In addition, we can take following form
\begin{align}
\phi(0,x) &= A_\phi [f_R(x-x_0) + f_R(x+x_0)], \\
\dot{\phi}(0,x) &= v A_\phi [\partial_x f_R (x+x_c) - \partial_x f_R (x-x_c),\\
\chi(0,x) &= A_\chi [ f_R(x-x_0)\cos\Delta\theta_{-} \nonumber \\
& + f_R(x+x_0)\cos\Delta\theta_{+} ], \\
\dot{\chi}(0,x) &= v A_\chi\big[\partial_x f_R(x+x_0)
                  \cos\Delta\theta_{-}\nonumber\\
 &\qquad-\partial_x f_R(x-x_0)\cos\Delta\theta_{+}\big],
\end{align}
which describes two bump collision.

If compact bumps live significantly longer than
matched Gaussian or periodic data, this would support the interpretation that dispersion and finite causal support suppress the nonlinear transfer channel.
If instead the lifetimes collapse when compared at fixed spectral centroid and amplitude, this would indicate that spectral content rather than compactness is the dominant control parameter. 

\subsection{Metastability observables and numerical limits}
\label{ssec:observables}

A central aim of the numerical analysis is to distinguish observed growth in the evolved fields from loss of convergence of the nonlinear algebraic solver. We therefore distinguish three quantities: the amplitude-threshold time $T_{\rm amp}$, the spectral-threshold time $T_{\rm spec}$, and the last converged terminal time $T_{\rm conv}$. Only the first two refer to operational changes in field diagnostics. The third refers to numerical continuation.

\subsubsection{Amplitude growth}

Define the normalized maximum field amplitude
\begin{equation}
\mathcal{A}(t)=
\frac{
\max_{x\in\Omega}
\left\{
|\phi(x,t)|,|\chi(x,t)|
\right\}
}{
\max_{x\in\Omega}
\left\{
|\phi(x,0)|,|\chi(x,0)|
\right\}
}.
\label{eq:normalized_amplitude}
\end{equation}

The amplitude-threshold time is
\begin{equation}
T_{\rm amp}
=
\inf\left\{
t:\mathcal{A}(t)>\eta_A
\right\},
\label{eq:Tamp}
\end{equation}
where $\eta_A>1$ specifies the adopted threshold.

In contrast with a mathematical blow-up time, $T_{\rm amp}$ records the first crossing of a finite diagnostic level. In oscillatory solutions, such a crossing can be transient.

The amplitude thresholds used in the parameter surveys are therefore stated explicitly in the associated figure captions.
We do not impose a universal $\eta_A$ across all initial-data families, and we do not interpret absolute threshold-crossing times obtained with different values of $\eta_A$ as directly comparable.

\subsubsection{Fourier-power diagnostics}

For periodic spatial data, define the field Fourier powers
\begin{equation}
P_\phi(k,t)=
|\widehat{\phi}(k,t)|^2,
\qquad
P_\chi(k,t)=
|\widehat{\chi}(k,t)|^2,
\label{eq:field_fourier_powers}
\end{equation}
and
\begin{equation}
P_{\rm tot}(k,t)=
P_\phi(k,t)+P_\chi(k,t).
\label{eq:total_field_power}
\end{equation}

These are squared Fourier amplitudes. They do not represent a decomposition of the signed Hamiltonian.

For a physical wavenumber cutoff $k_c$, we consider
\begin{equation}
R_{\rm IR/UV}(t)
=
\frac{
\displaystyle
\sum_{0<|k|\le k_c}
P_{\rm tot}(k,t)
}{
\displaystyle
\sum_{|k|>k_c}
P_{\rm tot}(k,t)
}.
\label{eq:RIRUV}
\end{equation}

The zero mode is excluded and can be monitored separately:
\begin{equation}
R_0(t)=
\frac{
P_{\rm tot}(0,t)
}{
\displaystyle\sum_{k\ne0}
P_{\rm tot}(k,t)
}.
\label{eq:R0}
\end{equation}

For cases in which the initial IR/UV ratio is meaningful, we define
\begin{equation}
\mathcal{G}_{\rm IR}(t)
=
\frac{
R_{\rm IR/UV}(t)+\epsilon_R
}{
R_{\rm IR/UV}(0)+\epsilon_R
},
\label{eq:GIR}
\end{equation}
where $\epsilon_R=10^{-12}$ is a regularization parameter.

The spectral-threshold time is
\begin{equation}
T_{\rm spec}
=
\inf\left\{
t:\mathcal{G}_{\rm IR}(t)>\Gamma
\right\}.
\label{eq:Tspec}
\end{equation}

The threshold $\Gamma$ and the Fourier cutoff are specified for each applicable survey. Typical values used in the present results are
$\Gamma=10$ and $\Gamma=20$.

A threshold crossing in $\mathcal{G}_{\rm IR}$ identifies a change in the relative distribution of field Fourier power. It is not sufficient by itself to establish a physical runaway or an inverse energy cascade.

\subsubsection{Selection of spectral cutoffs}

The periodic Fourier wavenumbers are indexed by
\begin{equation}
k_n=\frac{2\pi n}{L}.
\end{equation}

When a cutoff is specified by an integer Fourier mode $n_c$, the corresponding physical wavenumber is
\begin{equation}
k_c=\frac{2\pi n_c}{L}.
\label{eq:mode_cutoff_conversion}
\end{equation}

This distinction is essential because a dimensionless mode index and a physical wavenumber are not interchangeable.

For Gaussian initial data, the characteristic inverse-width prescription is
\begin{equation}
k_c=\frac{1}{4\ell}.
\label{eq:gaussian_cutoff}
\end{equation}

For colored-noise data, a possible spectral reference scale is
\begin{equation}
k_c=\sqrt{k_{n_1}k_{n_2}},
\end{equation}
although the actual Fourier-bin selection used in the post-processing must be reported consistently with the implemented cutoff.

For monochromatic plane waves, ordinary harmonic generation and the small number of initially populated modes can make the spectral threshold unsuitable as a general growth classifier.
We therefore use the amplitude criterion for the plane-wave and phase-correlated maps.

For broadband Gaussian, colored-noise, and compact-bump surveys, the spectral diagnostic provides a complementary measure of departure from the initial spectral distribution.

\subsubsection{Operational lifetime and censoring}

When both amplitude and spectral criteria are applicable, we define the operational lifetime
\begin{equation}
T_{\rm life}
=
\min\left(
T_{\rm amp},
T_{\rm spec}
\right),
\label{eq:Tlife}
\end{equation}
provided the applicable threshold is crossed during a converged
numerical evolution.

For surveys in which only the amplitude criterion is used, we set
\begin{equation}
T_{\rm life}=T_{\rm amp}.
\end{equation}

A threshold-crossing time is reported as measured only when the corresponding numerical solution is available and converged at that time.

If the run reaches its prescribed final time $T_{\rm final}$ without crossing the applicable growth criterion, the physical onset time is right-censored:
\begin{equation}
T_{\rm life}>T_{\rm final}.
\label{eq:right_censored}
\end{equation}

If the nonlinear solver fails before the threshold is crossed, the run is classified as solver-limited. We report its last converged time separately rather than assigning that time to a physical instability.

In particular,
\begin{equation}
T_{\rm conv}
\ne
T_{\rm life}
\end{equation}
in general.

The continuation time depends on the algorithm, nonlinear initial guess, discretization, and convergence tolerances, whereas the operational lifetime is defined from the evolved field diagnostics.

\subsubsection{Hamiltonian and numerical diagnostics}

To monitor the quality of the numerical evolution, we compute the relative total-Hamiltonian drift
\begin{equation}
\epsilon_H(t)=
\frac{|H(t)-H(0)|}
{1+|H(0)|}.
\label{eq:H_drift}
\end{equation}
We additionally monitor the nonlinear residual, convergence status, and field-amplitude histories.

These quantities are particularly important because large compensating normal- and ghost-sector energy changes can coexist with a finite total Hamiltonian.
The field Fourier powers and sector Hamiltonians serve different purposes: the former describe spatial spectral organization, while the latter describe the signed energy decomposition of the model.

Neither a solver termination nor a single spectral or sector-energy excursion is treated as proof of a finite-time singularity.

\subsubsection{Localization diagnostics}

For the localized initial-data families, we consider the second-moment radius
\begin{equation}
\sigma_x^2(t)=
\frac{
\displaystyle
\int_\Omega
(x-x_c)^2
\left(\phi^2+\chi^2\right)
\,\mathrm{d}x
}{
\displaystyle
\int_\Omega
\left(\phi^2+\chi^2\right)
\,\mathrm{d}x
}.
\label{eq:sigmax}
\end{equation}

Here $x_c$ denotes the reference center used in the diagnostic. On a periodic domain, the definition is sensitive to the chosen coordinate origin when substantial field weight approaches or crosses the periodic boundary.

For approximately uniform spatial field weight on a domain of length $L$ centered at $x_c$,
\begin{equation}
\sigma_x=
\frac{L}{\sqrt{12}}.
\label{eq:uniform_radius}
\end{equation}
This value provides a useful reference for interpreting the localization histories in Sec.~\ref{sec:phi6}.

A second-moment radius alone does not uniquely identify a localized oscillon. In particular, a persistent localized core should be distinguished from an extended or dispersive state whose second moment approaches the domain-scale reference value.

\section{Metastability landscapes in the quartic
ghost-normal model}
\label{sec:quartic_landscapes}

We first examine the finite-time dynamics of the quartically coupled ghost-normal scalar system. Unless otherwise stated, the simulations use
$m_\phi=m_\chi=1$, $\gamma=-1$, and
\begin{equation}
V_{22}(\phi,\chi)=\lambda_{22}\phi^2\chi^2.
\end{equation}

The parameter surveys provide complementary views of the dynamics. The plane-wave and Gaussian families vary characteristic frequency and spatial width, while the colored-noise family varies the distribution of initial Fourier power at fixed field root-mean-square amplitude. The phase-correlated runs examine relative waveform alignment, and the localized and compact-bump configurations test whether strongly localized initial perturbations remain bounded for extended intervals.

A common qualitative pattern emerges: increasing initial amplitude tends to produce earlier threshold crossing, while increasing characteristic wavenumber or ultraviolet spectral bias tends to delay the observed growth. Nevertheless, the results are not described by a universal instability time. Some runs cross an operational amplitude or spectral threshold, others reach the prescribed terminal time without a crossing, and several terminate because the nonlinear solver loses convergence while the fields remain bounded. Throughout this section, we maintain the distinction between these outcomes.

\subsection{Plane-wave initial conditions:
amplitude--wavenumber landscape}
\label{ssec:quartic_plane}

\begin{figure}[!H]
\centering
\includegraphics[width=\linewidth]{
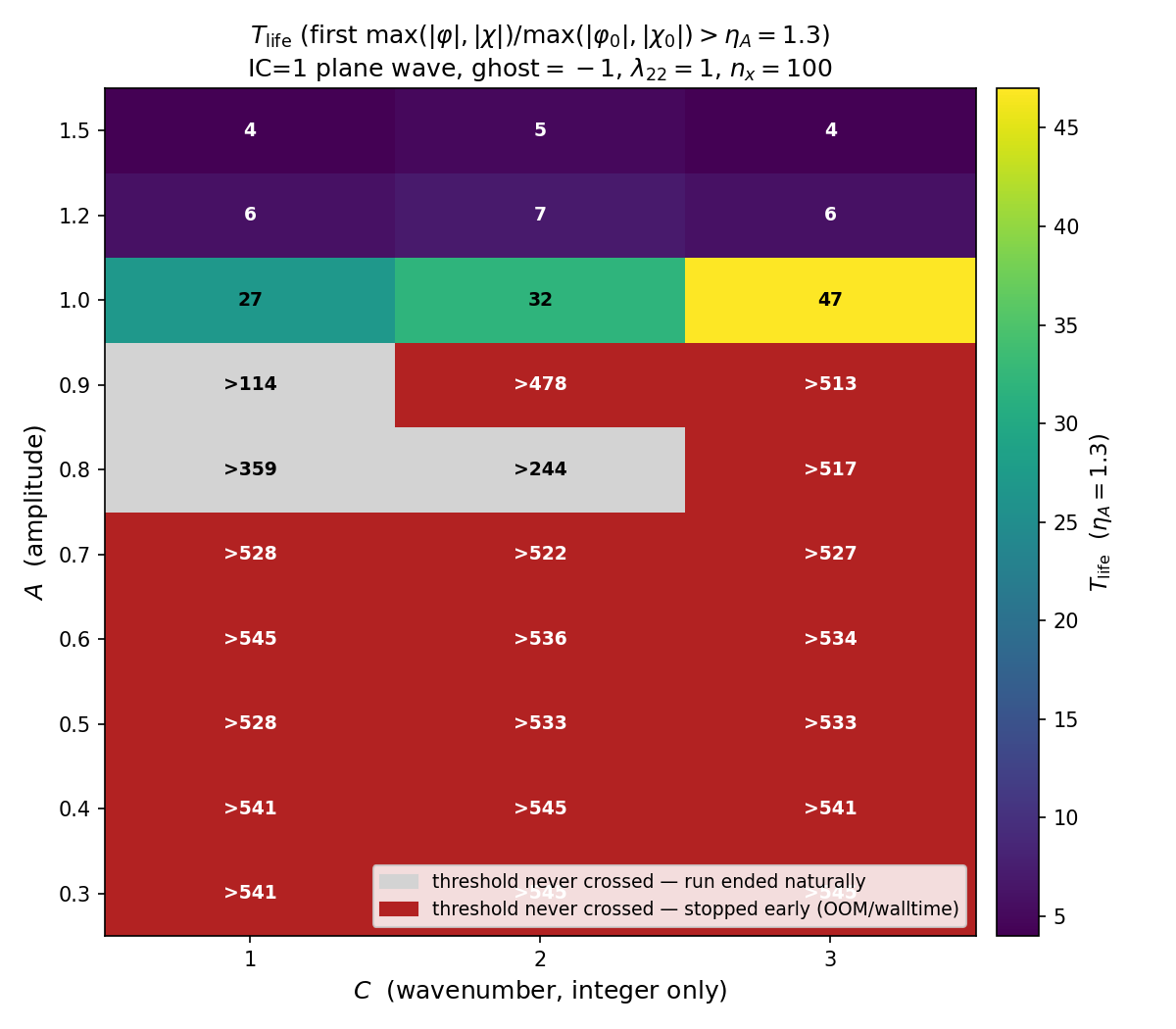
}
\caption{
Plane-wave amplitude--wavenumber map for the quartic ghost-normal model. The plotted onset time is $T_{\rm amp}$ with $\eta_A=1.3$. Red cells reached the externally imposed terminal-time cap without crossing the amplitude threshold and are right-censored.
Gray cells lost nonlinear-solver convergence before the threshold was crossed and are classified as solver-limited. Neither category is assigned a physical instability time.
}
\label{fig:IC_1_Heatmap}
\end{figure}

The plane-wave survey provides a controlled comparison of initial amplitude and carrier wavenumber. We vary the common amplitude $A$
and the dimensionless wavenumber parameter
\begin{equation}
C=\frac{kL}{2\pi},
\end{equation}
while retaining the relative wave configuration defined in Sec.~\ref{ssec:id_plane}.

For this initial-data family, we use the amplitude threshold $T_{\rm amp}$ rather than the spectral threshold $T_{\rm spec}$ as the primary onset diagnostic. A single-carrier configuration can generate higher harmonics during otherwise bounded nonlinear motion. Consequently, Fourier-power redistribution across a fixed cutoff need not coincide with the onset of strong field growth.

Figure~\ref{fig:IC_1_Heatmap} shows that large-amplitude configurations typically cross the adopted growth threshold earlier than smaller-amplitude configurations.
At comparable amplitude, increasing the carrier wavenumber generally delays threshold crossing over the measured portion of the parameter space.
The ordering is consistent with the frequency dependence previously reported for interacting ghost-normal scalar fields
by~\cite{Deffayet2025}.

An equally important observation concerns the distinction between a physical threshold event and numerical termination. Several configurations, including the parameter pairs
\begin{equation}
(A,C)=(0.9,1),\quad (0.8,1),\quad (0.8,2),
\end{equation}
terminate numerically while their maximum field amplitudes remain close to the initial values.
These runs do not provide evidence of runaway-like growth before the solver terminates. They demonstrate that the last converged terminal time $T_{\rm conv}$ cannot generally be identified with $T_{\rm amp}$.

The main conclusion from the plane-wave survey is therefore an ordering of operational growth times:
\begin{equation}
A\uparrow
\quad\Longrightarrow\quad
T_{\rm amp}\downarrow,
\qquad
k\uparrow
\quad\Longrightarrow\quad
T_{\rm amp}\uparrow,
\label{eq:plane_lifetime_ordering}
\end{equation}
within the regions where the threshold is reached and the ordering is resolved.

We do not interpret Eq.~\eqref{eq:plane_lifetime_ordering} as a universal scaling law or as proof that every initially bounded configuration eventually becomes singular.

\subsection{Gaussian packets: amplitude and spatial-width dependence}
\label{ssec:quartic_gaussian}

Gaussian packets extend the plane-wave comparison to localized initial profiles. The common initial amplitude $A$ is varied together with the width $\ell$, represented by
\begin{equation}
C=\frac{L}{8\pi\ell}.
\end{equation}
The simulation uses co-propagating initial data with both packet centers at $x_\phi=x_\chi=0.5L$.

Increasing $C$ corresponds to reducing the packet width and increasing its characteristic inverse-width spectral scale. This changes not only the initial Fourier distribution but also the spatial gradient contribution to the initial energy.
The resulting comparison should therefore be interpreted as a joint width--spectral dependence rather than a controlled
variation of only one spectral moment.

\begin{figure}[H]
\centering
\includegraphics[width=\linewidth]{
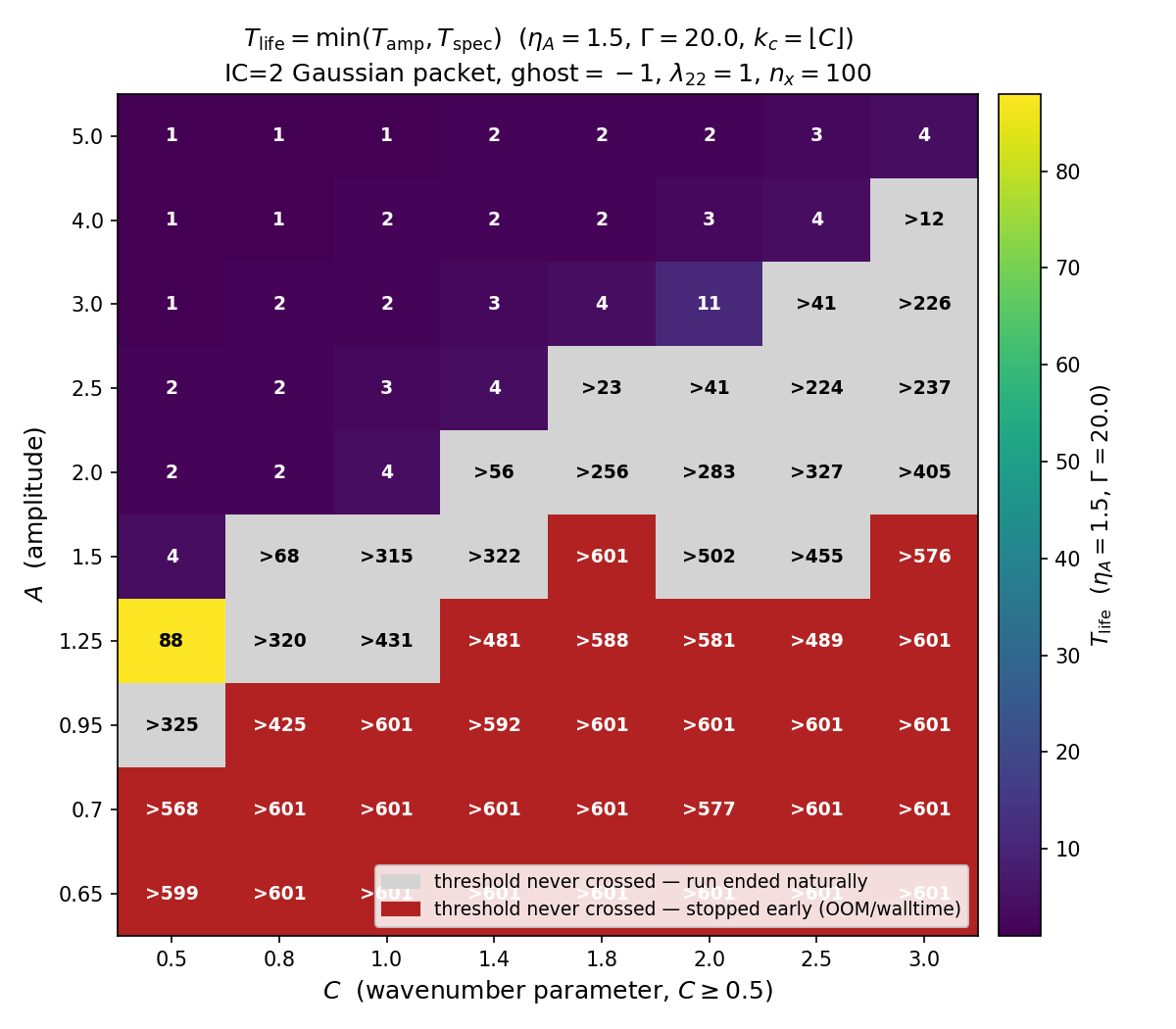
}
\caption{Operational threshold-crossing map for co-propagating Gaussian packets as a function of amplitude $A$ and dimensionless inverse-width parameter $C=k_{\rm char}L/(2\pi)$, with $k_{\rm char}=1/(4\ell)$. We use $\eta_A=1.5$ and $\Gamma=20$ and report $T_{\rm life}=\min(T_{\rm amp},T_{\rm spec})$ when an applicable threshold is reached during a converged evolution. Right-censored and solver-limited cases are distinguished in the map.}
\label{fig:IC_2_Heatmap}
\end{figure}

The map in Fig.~\ref{fig:IC_2_Heatmap} exhibits the same broad qualitative ordering seen in the plane-wave survey. At fixed width, larger amplitudes generally lead to earlier threshold crossings. At comparable amplitude, broader packets, corresponding to smaller $C$, tend to cross earlier than narrower packets.
This is consistent with the spectral organization of the initial data. Broader Gaussian profiles have greater relative weight at lower wavenumbers, whereas narrower profiles contain a wider range
of spatial frequencies.

However, because the initial width also changes the gradient energy and degree of localization, the heatmap alone does not establish that the spectrum is the unique cause of the observed ordering.
For Gaussian data, the spectral ratio $R_{\rm IR/UV}$ provides more information than for an initially monochromatic wave. Nevertheless, strongly oscillatory solutions may temporarily cross either the amplitude or spectral threshold without subsequently exhibiting runaway.

Accordingly, the plotted quantities are operational onset indicators, not mathematically defined finite-time singularity times. The Gaussian results establish a correlation between the initial amplitude--width configuration and the time at which the evolution departs from the adopted bounded regime.

Section~\ref{sec:spectral_mechanism} examines time-resolved Fourier spectra to determine whether low-wavenumber enhancement can precede the amplitude growth in a representative pair of
Gaussian configurations.

\subsection{Colored-noise spectra: spectral tilt and phase variability}
\label{ssec:quartic_colored}

The colored-noise result provides a more direct investigation of initial spectral bias. The field RMS amplitude is fixed by the normalization of Sec.~\ref{ssec:id_colored}, while the spectral tilt $n_s$ redistributes the initial field Fourier power.
Negative spectral tilts place greater relative power at low wavenumbers, whereas positive tilts favor higher wavenumbers. For these simulations, we use $s_n=-1$ and $n_2=32$ and vary both amplitude and spectral tilt.

\begin{figure}[H]
\centering
\includegraphics[width=\linewidth]{
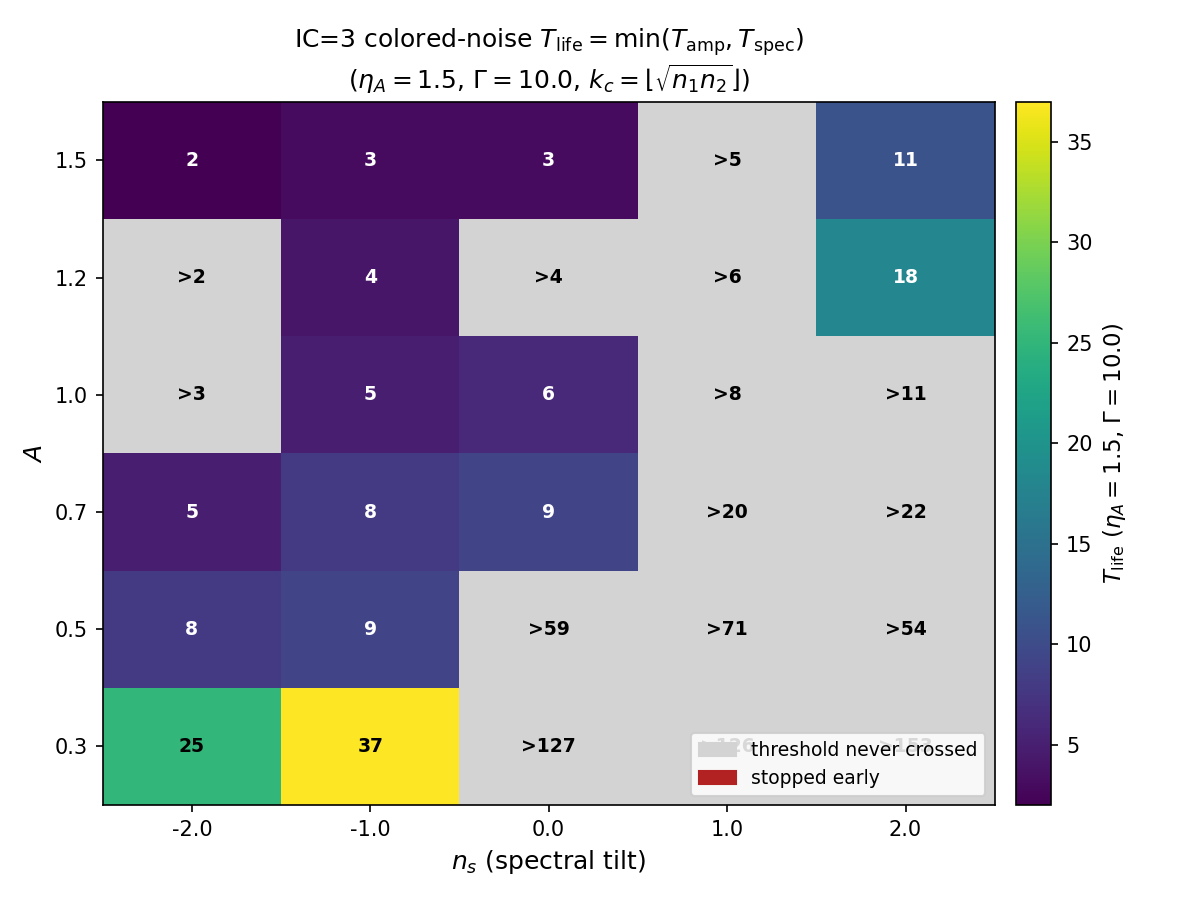
}
\caption{Representative colored-noise threshold-crossing map over initial RMS amplitude $A$ and spectral tilt $n_s$. The diagnostics use $\eta_A=1.5$ and $\Gamma=10$. Negative $n_s$ is infrared biased, whereas positive $n_s$ is ultraviolet biased. The displayed map corresponds to one random-phase realization; two additional realizations were used to assess the qualitative ordering and the variability of onset times.}
\label{fig:IC_3_heatmap_Tlife_eta1.5_Gamma10.0}
\end{figure}

Figure~\ref{fig:IC_3_heatmap_Tlife_eta1.5_Gamma10.0} shows a clear qualitative dependence on the initial spectral tilt. At comparable amplitude, infrared-biased configurations tend to cross the adopted growth criteria earlier than ultraviolet-biased configurations.
The selected time-series diagnostics support this interpretation. In the infrared-biased cases, the field-amplitude envelope typically exhibits more pronounced growth. By contrast, representative strongly UV-biased cases can show initial amplitude reduction followed by approximately bounded evolution over the resolved interval.
The spectral ordering persists qualitatively across three independent random-phase realizations. However, the onset times are not independent of the specific phase realization.

Table~\ref{tab:ic3-seed-ensemble} summarizes representative parameter choices for which the run-by-run onset information is available.
\begin{table*}[!t]
\centering
\caption{Three seed colored-noise threshold-crossing summary. The diagnostic uses $\eta_A=1.5$ and $\Gamma=10$, with the reported Fourier cutoff corresponding to the existing analysis convention. Entries of the form $>t$ denote runs in which the applicable threshold was not crossed before the available terminal time. These entries are lower bounds, not measured onset times. Median and range are reported only when all three onset times are measured.}
\label{tab:ic3-seed-ensemble}
\begin{tabular}{cccccccc}
\toprule
$A$ & $n_s$ &
Seed 1 & Seed 2 & Seed 3 &
Measured & Median & Range\\
\midrule
$0.3$ & $-2$ &
$70$ & $255$ & $15$ &
$3/3$ & $70$ & $[15,255]$\\
$0.3$ & $-1$ &
$31$ & $272$ & $60$ &
$3/3$ & $60$ & $[31,272]$\\
$0.3$ & $+2$ &
$>147$ & $>181$ & $>147$ &
$0/3$ & --- & ---\\
$0.5$ & $-2$ &
$13$ & $8$ & $7$ &
$3/3$ & $8$ & $[7,13]$\\
$0.5$ & $-1$ &
$18$ & $20$ & $17$ &
$3/3$ & $18$ & $[17,20]$\\
$0.5$ & $+2$ &
$>52$ & $>67$ & $>58$ &
$0/3$ & --- & ---\\
\bottomrule
\end{tabular}
\end{table*}
At $A=0.3$ and $n_s=-2$, the three measured threshold times are $15$, $70$, and $255$. For $n_s=-1$ at the same amplitude, the corresponding times are $31$, $60$, and $272$.
Thus, the measured onset times can vary considerably between realizations even when the spectral tilt and RMS field amplitude are held fixed. 
At $A=0.5$, the spread is smaller for the tabulated infrared-biased samples: the values are $7$--$13$ for $n_s=-2$ and $17$--$20$ for $n_s=-1$.

The displayed ultraviolet-biased cases with $n_s=+2$ do not cross the selected thresholds before their available terminal times. Consequently, their exact onset times and corresponding medians are not determined. These results support two distinct conclusions.
First, the initial spectral tilt organizes the qualitative finite-time response: UV-biased initial fields tend to remain below the adopted growth thresholds for longer periods.
Second, random phase realizations can substantially affect the quantitative threshold-crossing time.
The latter observation is not inconsistent with the former. Rather, it demonstrates that the spectral power distribution and the detailed phase structure both contribute to the finite-time nonlinear dynamics.

Finally, fixing the field RMS does not fix the initial gradient, kinetic, or total Hamiltonian contributions. The colored-noise survey therefore establishes a comparison at fixed RMS and varying spectral tilt, not a comparison at fixed total energy.

\subsection{Phase-correlated two-field waves}
\label{ssec:quartic_phase}
We next examine how relative waveform alignment affects finite-time evolution. The initial data of Sec.~\ref{ssec:id_phase} contain a common carrier wavenumber and prescribed normal-field amplitude $A$, while varying the ghost-to-normal amplitude ratio $r$ and relative phase $\Delta\phi$.

We fix
\begin{equation}
A=1,
\qquad
C=\frac{kL}{2\pi}=1,
\end{equation}
and consider co-propagating and counter-propagating configurations separately.

As in the plane-wave survey, we use the amplitude threshold rather than the spectral threshold to characterize growth, since ordinary harmonic generation can lead to spectral threshold crossings without corresponding strong field-amplitude growth.

\begin{figure}[H]
\centering
\includegraphics[width=0.9\linewidth]{
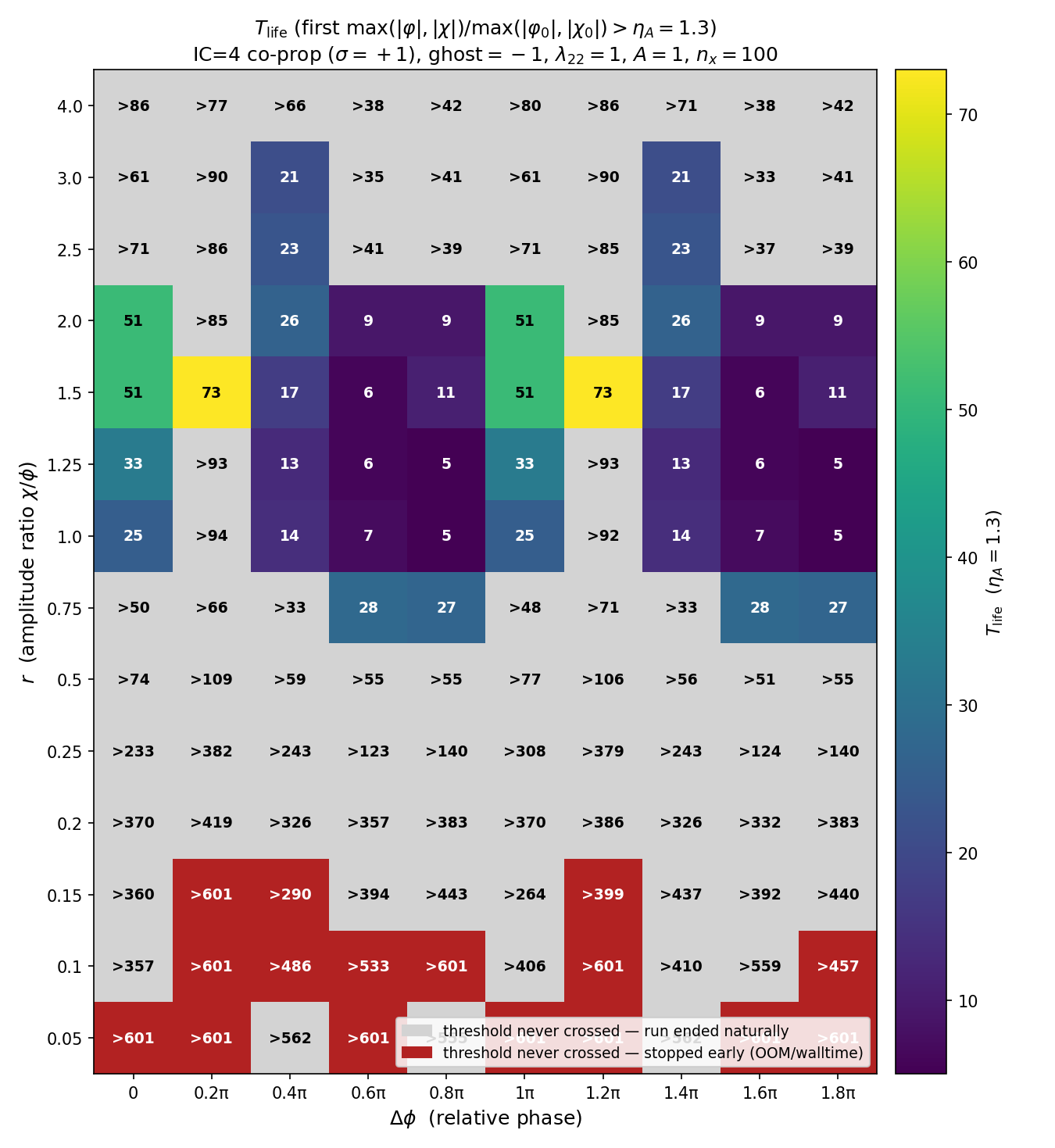
}
\caption{Co-propagating phase-correlated plane-wave survey. The operational amplitude threshold uses $\eta_A=1.3$. Several configurations terminate because the nonlinear solver loses convergence before the amplitude threshold is reached.}
\label{fig:IC_4_co-prop_heatmap}
\end{figure}

\begin{figure}[H]
\centering
\includegraphics[width=\linewidth]{
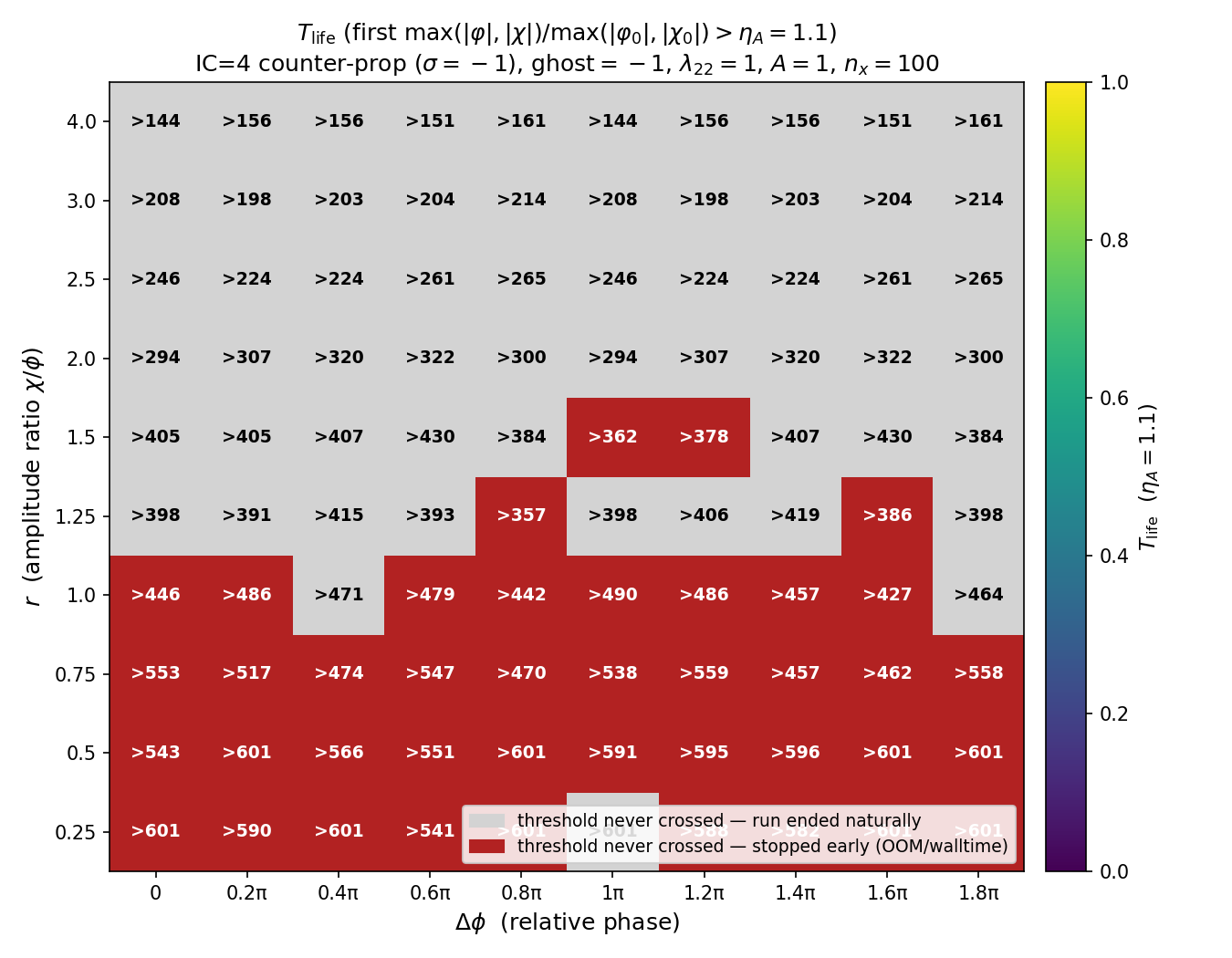
}
\caption{Counter-propagating phase-correlated plane-wave survey. The plotted amplitude criterion uses $\eta_A=1.1$. Most configurations do not cross this sensitive threshold before the imposed stopping time or nonlinear-solver termination. Because this threshold differs from the co-propagating survey, the absolute displayed onset times are not directly comparable between the maps.}
\label{fig:IC_4_counter-prop_heatmap}
\end{figure}

Figures~\ref{fig:IC_4_co-prop_heatmap} and~\ref{fig:IC_4_counter-prop_heatmap} suggest a qualitative difference between the two propagation configurations.
The co-propagating map shows more visible dependence on relative phase, while the counter-propagating map exhibits predominantly low-amplitude oscillatory behavior over the resolved intervals.
Increasing $r$ also tends to reduce the available converged evolution time over parts of the parameter space.

A possible physical interpretation is that co-propagating field patterns maintain their relative alignment over longer intervals, allowing their nonlinear interaction to remain sensitive to the initial phase.
Counter-propagating patterns instead change their relative spatial alignment during evolution, which may reduce the persistence of phase-sensitive effects.
However, these observations should not be interpreted as proof of a resonance or a measured phase-dependent instability boundary.

Many configurations are solver-limited before strong growth is observed. Furthermore, the maps use different amplitude thresholds, preventing a direct quantitative comparison of their absolute threshold-crossing times.
The supported conclusion is therefore that the resolved finite-time dynamics displays qualitative sensitivity to initial phase structure, especially in the co-propagating family.

\subsection{Localized time-symmetric seeds: quartic-model controls}
\label{ssec:quartic_localized}
The localized time-symmetric $\sech$ data provide a useful control family for examining the roles of the ghost kinetic sign and nonlinear coupling. We compare configurations with $\gamma=\pm1$ and $\lambda_{22}=0,1$ while keeping the initial profile fixed. The simulations use $r=1$, $\Delta\phi=0$, and no additional carrier, $k_0=0$.

\begin{table}[!htbp]
\centering
\caption{Last converged terminal times for localized time-symmetric quartic-model controls, with $r=1$, $\Delta\phi=0$, and $k_0=0$. The values are numerical continuation times and are not automatically physical instability times. The computational cap is 601.}
\label{tab:IC_5-lam-v-gam}
\begin{tabular}{ccc}
\toprule
& $\gamma=+1$ & $\gamma=-1$\\
\midrule
$\lambda_{22}=0$ & 601 & 601\\
$\lambda_{22}=1$ & 315 & 430\\
\bottomrule
\end{tabular}
\end{table}

\begin{figure}[H]
\centering
\includegraphics[width=\linewidth]{
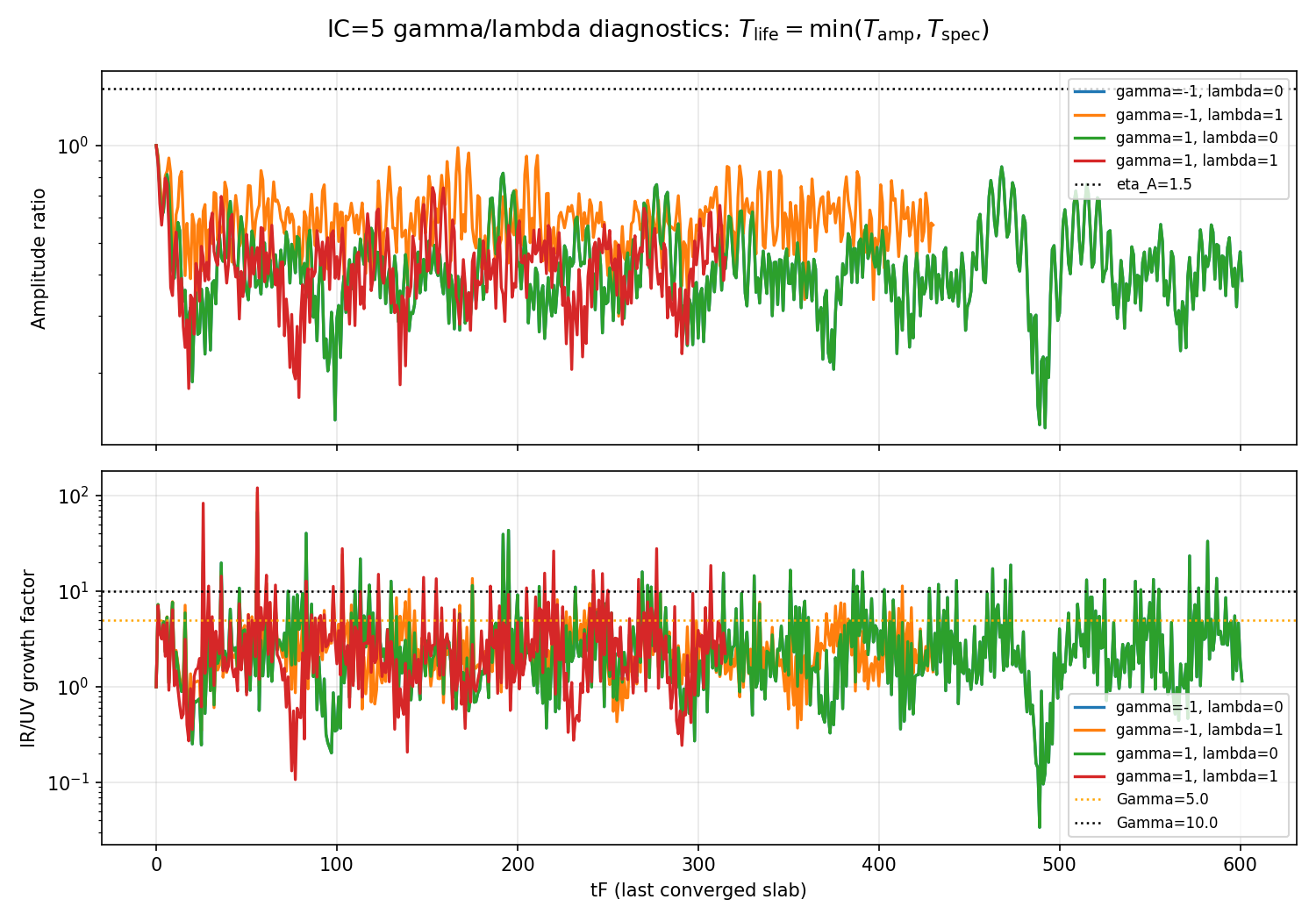
}
\caption{Amplitude and IR-to-UV Fourier-power diagnostics for the quartic localized controls. The nonlinear runs terminating at finite $T_{\rm conv}$ do not show corresponding runaway-like amplitude growth or pronounced spectral-threshold behavior over their available resolved intervals. Their termination is therefore classified as solver-limited rather than as measured physical blow-up.}
\label{fig:IC_5_lam-v-gam}
\end{figure}

Table~\ref{tab:IC_5-lam-v-gam} shows that the nonlinear $\lambda_{22}=1$ cases terminate at smaller $T_{\rm conv}$ than the corresponding uncoupled controls. However, the diagnostic histories in Fig.~\ref{fig:IC_5_lam-v-gam} do not show the strong field-amplitude growth that would justify identifying these terminal times with physical runaway.

In particular, the $\gamma=+1$ control terminates earlier than the corresponding $\gamma=-1$ run in this comparison.
This does not establish that the positive-energy theory is physically less stable than the ghost-normal theory. It instead highlights the sensitivity of nonlinear-solver convergence to the algebraic problem and its initial iterate.
The quartic localized controls thus reinforce a central methodological conclusion of the paper: solver termination alone cannot be used as an instability classifier.
They also provide a reference for the higher-order interaction
model examined in Sec.~\ref{sec:phi6}.

\subsection{Smooth compactly supported perturbations}
\label{ssec:quartic_compact}
The final quartic-model initial-data family consists of smooth compactly supported, time-symmetric perturbations.
We use the single-bump data introduced in Sec.~\ref{IC:compact-bumps}, with
\begin{equation}
x_0=0.5L,
\qquad
R=0.1L,
\end{equation}
and initially vanishing time derivatives.

The survey varies the amplitude parameter and the factor $\cos\Delta\theta$ multiplying the initial ghost-field profile.

\begin{figure}[H]
\centering
\includegraphics[width=\linewidth]{
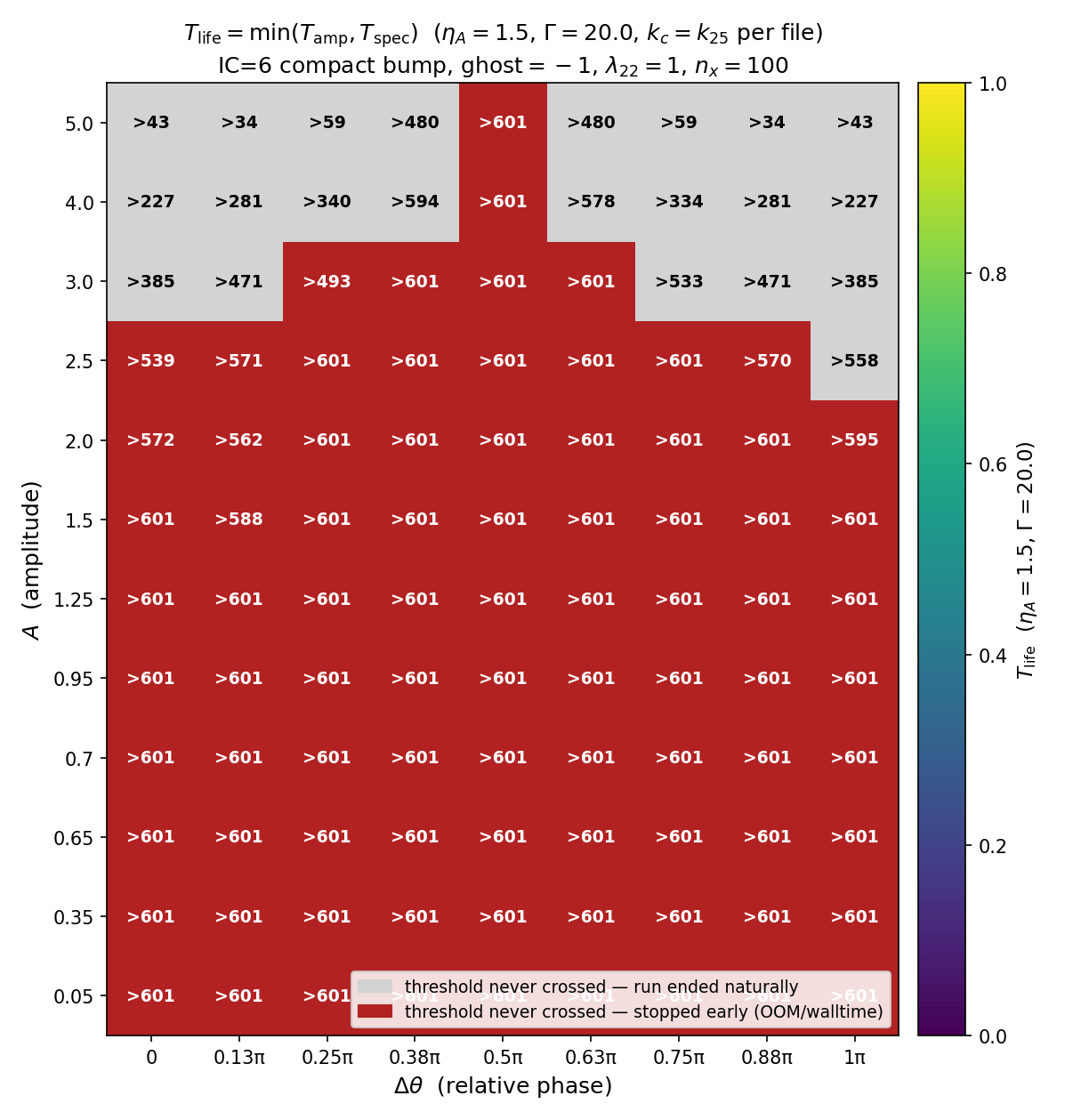
}
\caption{Operational threshold-crossing map for initially compact single-bump perturbations with $\dot{\phi}=\dot{\chi}=0$, $x_0=0.5L$, and $R=0.1L$. The adopted thresholds are $\eta_A=1.5$ and $\Gamma=20$. A substantial portion of the small- and moderate-amplitude parameter space reaches the imposed terminal-time cap without crossing either applicable growth threshold.}
\label{fig:IC_6_Heatmap}
\end{figure}

Figure~\ref{fig:IC_6_Heatmap} exhibits a broad region of long resolved evolution without amplitude or spectral threshold crossing at small and moderate initial amplitudes.
Only at larger amplitudes does a substantial part of the map move away from this long-lived region.

This result is qualitatively consistent with the importance of small initial data in recent mathematical work on selected ghost-normal wave systems~\cite{Held2025b}. The comparison must nevertheless remain carefully qualified.
Our simulations are performed over finite times on a periodic computational domain. They do not establish global-in-time boundedness or scattering to infinity for compact data.

Furthermore, the apparent longer-lived behavior near $\Delta\theta=\pi/2$ is not evidence for a separate phase-stabilization mechanism.

The initial-data ansatz contains
\begin{equation}
\chi(x,0)\propto\cos\Delta\theta,
\end{equation}
so this value suppresses the initially excited ghost field.

For the quartic interaction, setting the ghost field and its initial time derivative to zero also eliminates the nonlinear source that would otherwise excite it from the normal sector.
The strongest conclusion from the compact-bump survey is therefore the broad range of initially localized, low-amplitude configurations that remain below the adopted growth thresholds over the simulated interval.
The result complements the plane-wave, Gaussian, and colored-noise surveys by showing that long-lived bounded numerical behavior also occurs for spatially compact initial perturbations.

\section{Spectral organization of finite-time metastability}
\label{sec:spectral_mechanism}
The results in Sec.~\ref{sec:quartic_landscapes} show that the initial spectral distribution influences the finite-time response of the quartic ghost-normal system.
In the plane-wave survey, increasing the carrier wavenumber generally delays the adopted amplitude-growth threshold. Narrower Gaussian packets also tend to remain below the selected growth criteria for longer intervals, while UV-biased colored-noise initial data generally exhibit longer resolved bounded evolution than IR-biased configurations.
These observations establish a qualitative relationship between the initial spectral organization and the finite-time dynamics. They do not, by themselves, identify a specific nonlinear transfer mechanism. 

To obtain a more direct dynamical comparison, we examine the time-resolved field Fourier power and amplitude histories for two Gaussian configurations from the survey in Sec.~\ref{ssec:quartic_gaussian}. Both runs use the same characteristic inverse-width parameter $C=2$, but have different initial amplitudes:
\begin{equation}
(A,C)=(3,2)
\qquad\text{and}\qquad
(A,C)=(0.65,2).
\end{equation}
The first configuration exhibits rapid growth over its resolved evolution, whereas the second remains comparatively bounded over the simulated interval.

\begin{figure*}[!t]
\centering
\includegraphics[width=0.92\textwidth]{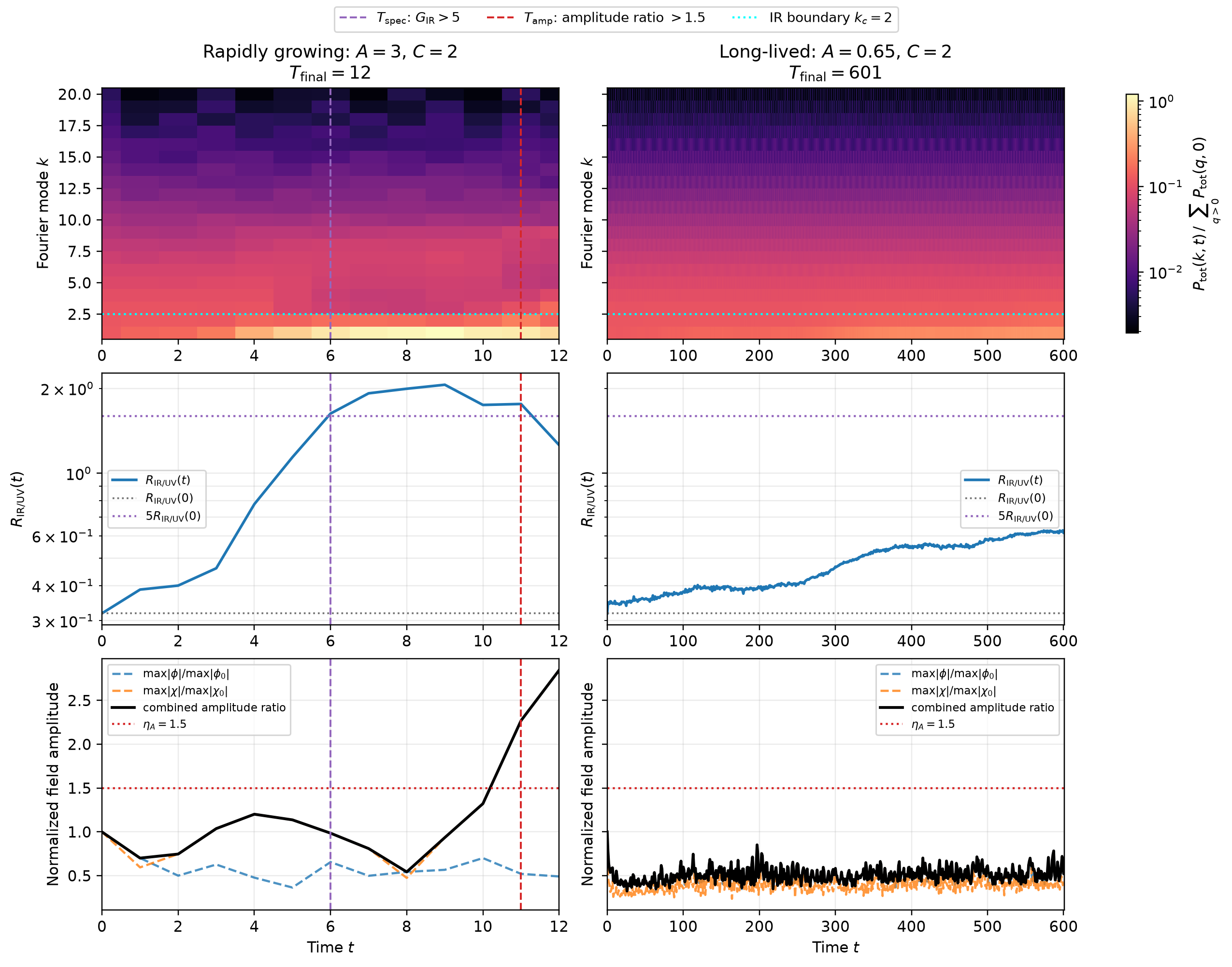}
\caption{Comparison of a rapidly growing Gaussian evolution (left, $A=3$, $C=2$) and a long-lived evolution (right, $A=0.65$, $C=2$). The upper panels show the time-dependent field Fourier-power distribution, the middle panels show the infrared-to-ultraviolet power ratio, and the lower panels show normalized field amplitudes. For the rapidly growing configuration, the selected spectral threshold is crossed near $t\simeq6$, before the amplitude threshold near $t\simeq11$.
The long-lived comparison exhibits substantially smaller low-wavenumber enhancement and no corresponding strong field-amplitude excursion over the available interval. The plotted Fourier powers are not signed Hamiltonian-energy densities.}
\label{fig:gaussian_spectral_precursor}
\end{figure*}

Figure~\ref{fig:gaussian_spectral_precursor} shows an important distinction between the two trajectories.
In the rapidly growing configuration, the IR-to-UV Fourier-power ratio increases before the substantial field-amplitude excursion. Using the diagnostic levels displayed in the figure, the spectral threshold is crossed at approximately $t=6$, whereas the amplitude threshold is reached at approximately $t=11$.
The temporal ordering suggests that redistribution of field power toward longer wavelengths can precede the onset of rapid amplitude growth.

In contrast, the long-lived comparison does not exhibit a similarly pronounced increase in the infrared power ratio. Its field amplitudes also remain comparatively close to the initial scale over the resolved interval.
The difference is consistent with the broader lifetime maps: configurations that develop substantial low-wavenumber enhancement can enter a rapidly growing regime, whereas those without comparable enhancement may remain bounded for longer times.

The comparison therefore supplies a concrete example of a \emph{spectral precursor} to strong field growth.
Nevertheless, several qualifications are essential.
First, the plotted diagnostic is the ratio of squared field Fourier amplitudes,
\begin{equation}
P_{\rm tot}(k,t)
=
|\widehat{\phi}(k,t)|^2
+
|\widehat{\chi}(k,t)|^2.
\end{equation}
It is not a decomposition of the signed Hamiltonian into spectral
energy contributions.
Consequently, enhancement of low-wavenumber Fourier power should not automatically be identified with a transfer of positive or negative Hamiltonian energy between Fourier modes.

Second, the time ordering observed in Fig.~\ref{fig:gaussian_spectral_precursor} does not establish that infrared enhancement is sufficient to cause runaway, nor that it is a necessary precursor in every initial-data family.

Third, the two configurations have different initial amplitudes and therefore also differ in their initial nonlinear interaction strength.
The comparison demonstrates a dynamical association between spectral reorganization and field growth, rather than isolating a unique causal mechanism.
A complete identification of the underlying transfer process would require mode-resolved signed-energy budgets or nonlinear spectral flux diagnostics. Such an analysis remains an interesting direction for future work.

In summary, the present results  support two conclusions.
The initial spectral distribution is an important organizing parameter for finite-time behavior, and significant low-wavenumber power enhancement can precede rapid field growth in the quartic ghost-normal system.
Together, these findings provide a physical interpretation of the metastability landscapes without requiring the assumption of a universal inverse cascade.

\section{Lifted $\phi^6$ model: potential-dependent continuation dynamics}
\label{sec:phi6}

We next investigate how the finite-time response changes when the quartic interaction is replaced by a different polynomial potential.
The motivation is provided by ordinary scalar theories with nonlinear potentials that can support long-lived localized oscillatory configurations, including oscillons
\citep{GleiserSicilia2009,Amin2010,
HindmarshSalmi2008}.

Our aim is to examine whether the corresponding lifted interaction changes the accessible finite-time behavior in a coupled ghost-normal theory.
We emphasize from the outset that a long converged simulation does not by itself establish the existence of a localized oscillon or a physically stable solution.

\subsection{Interaction and linearized dynamics}
\label{ssec:phi6_model}
We consider
\begin{equation}
U(W)
=
\frac12m_6^2 W
-\frac{\lambda}{4}W^2
+\frac{g}{6}W^3,
\qquad
W=\phi^2+\chi^2.
\label{eq:phi6_potential}
\end{equation}

In the numerical implementation, the quadratic term in Eq.~\eqref{eq:phi6_potential} is added to the bare mass terms already present in the equations of motion.

The resulting field equations are
\begin{align}
\Box\phi
+
\left(
m_\phi^2+m_6^2
-\lambda W+gW^2
\right)\phi
&=0,
\label{eq:phi6_eom_phi}
\\
\Box\chi
+
\left[
m_\chi^2+
\gamma\left(
m_6^2-\lambda W+gW^2
\right)
\right]\chi
&=0.
\label{eq:phi6_eom_chi}
\end{align}

For the parameters
\begin{equation}
m_\phi=m_\chi=m_6=\lambda=g=1,
\qquad
\gamma=-1,
\end{equation}
these reduce to
\begin{align}
\Box\phi+(2-W+W^2)\phi&=0,
\label{eq:phi6_eom_phi_numerical}
\\
\Box\chi+(W-W^2)\chi&=0.
\label{eq:phi6_eom_chi_numerical}
\end{align}

The linearized equations around the trivial vacuum are therefore
\begin{align}
\Box\phi+2\phi&=0,
\\
\Box\chi&=0.
\end{align}

In particular, the normal field has squared mass $m_{\phi,\mathrm{eff}}^2=2$, while the ghost field is massless at the linearized level: $m_{\chi,\mathrm{eff}}^2=0$.
This distinction is important when comparing the lifted potential with the quartic baseline, for which $m_\phi=m_\chi=1$.

The two systems differ both in their higher-order nonlinear couplings and in their linear dispersion relations.
Accordingly, their finite-time behavior must be interpreted as a comparison between complete field-theory models, not as an isolated measurement of sextic-interaction stabilization.

\subsection{Reference amplitude and initial data}
\label{ssec:phi6_reference}
The ordinary single-field potential motivating Eq.~\eqref{eq:phi6_potential} contains an attractive quartic contribution and a stabilizing higher-order term.
One useful reference condition is that the nonlinear restoring force be softer than its small-amplitude quadratic contribution over some range of field values.

For the polynomial interaction, the derivative is
\begin{equation}
\frac{\mathrm{d}U}{\mathrm{d}W}
=
\frac12
\left(
m_6^2-\lambda W+gW^2
\right).
\end{equation}

The nonlinear correction to the quadratic restoring coefficient is negative when
\begin{equation}
-\lambda W+gW^2<0,
\end{equation}
which, for positive
$\lambda$ and $g$, gives
\begin{equation}
0<W<\frac{\lambda}{g}.
\label{eq:phi6_softening}
\end{equation}

For the time-symmetric localized initial data of Eq.~\eqref{IC: oscillon-like}, the central field values give
\begin{equation}
W_{\rm max}
=
A^2
\left(
1+r^2\cos^2\Delta\phi
\right).
\end{equation}

This motivates the reference amplitude
\begin{equation}
A_c
=
\sqrt{
\frac{\lambda}{
g(1+r^2\cos^2\Delta\phi)
}
}.
\label{eq:phi6_reference_amplitude}
\end{equation}

For the simulations below, $r=1$, $\Delta\phi=0$, and $\lambda=g=1$, so that
\begin{equation}
A_c=\frac{1}{\sqrt{2}}.
\end{equation}

We emphasize that $A_c$ is a reference scale associated with the elementary softening condition in the potential.
It is not a mathematical existence condition for an oscillon in the coupled ghost-normal theory. Nor is it a demonstrated boundary between stable and unstable initial data.
The parameter scans therefore use $A/A_c$ simply as a convenient dimensionless amplitude coordinate.

\subsection{Amplitude--width continuation landscape}
\label{ssec:phi6_landscape}
We use time-symmetric localized initial data with no carrier, $k_0=0$, and vary the initial width $\sigma$ and normalized amplitude $A/A_c$.

The maximum terminal time in the numerical continuation sequence is
\begin{equation}
T_{\rm max}=601.
\end{equation}

\begin{figure*}[!t]
\centering
\includegraphics[width=0.60\textwidth]{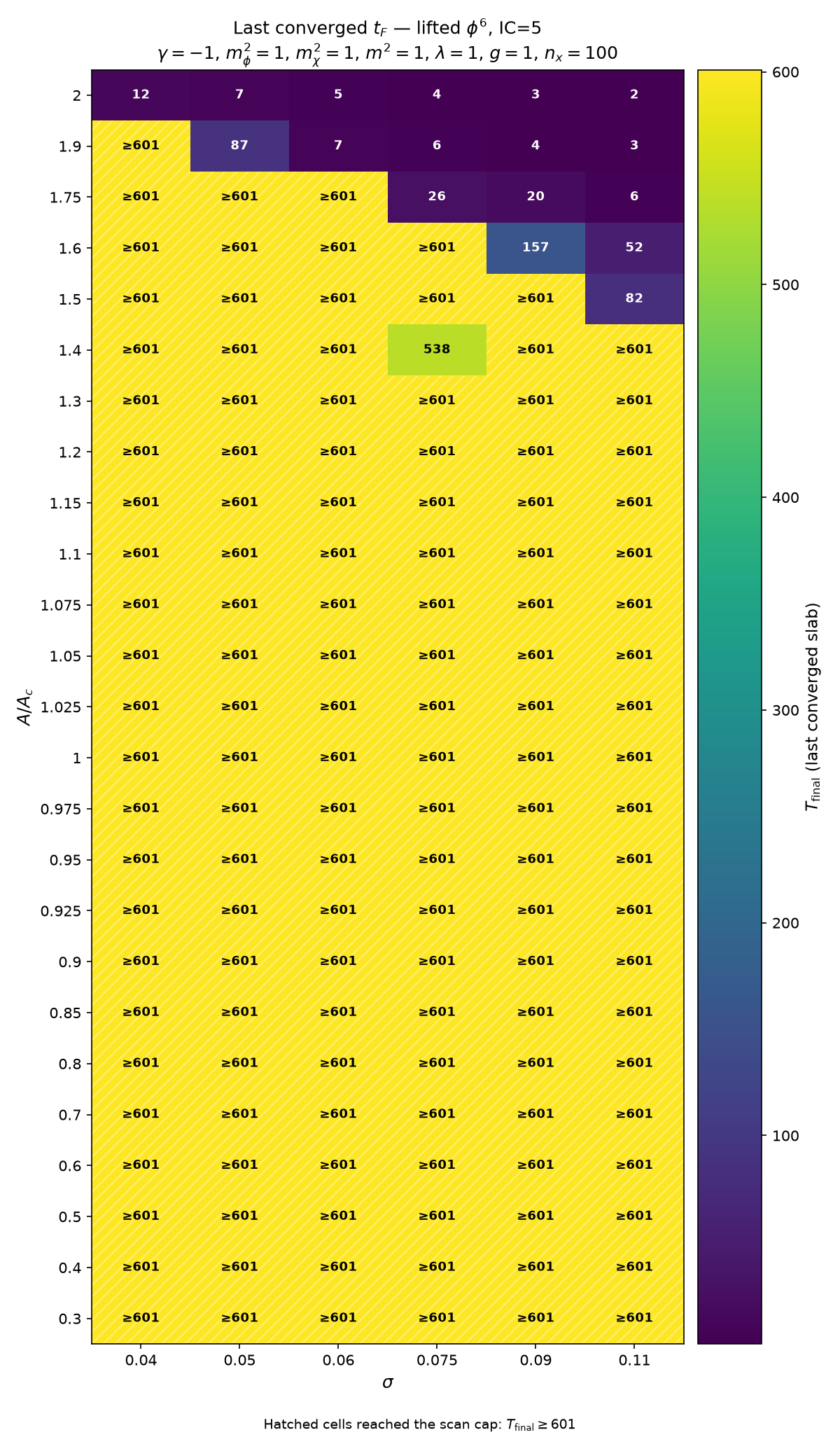}
\caption{Last converged terminal time $T_{\rm conv}$ for the lifted $\phi^6$ model as a function of initial width $\sigma$ and normalized amplitude $A/A_c$. The simulations use $m_\phi=m_\chi=m_6=\lambda=g=1$, $\gamma=-1$, $r=1$, and $\Delta\phi=0$. The map is obtained with warm-start continuation. Cells reaching the imposed cap $T_{\rm max}=601$ have no measured physical instability time. The figure therefore characterizes the accessible numerical continuation landscape rather than a physical lifetime boundary.}
\label{fig:phi6_continuation_map}

\end{figure*}

Figure~\ref{fig:phi6_continuation_map} shows a broad region of parameter space whose numerical solutions remain converged up to the prescribed terminal-time cap.
This region includes initial amplitudes substantially larger than the reference scale $A_c$.
At larger amplitudes, the last converged terminal time can decrease sharply, with some dependence on the initial width.

The resulting structure is qualitatively different from a simple transition at $A=A_c$.
Indeed, the long-converged region extends beyond that reference amplitude, demonstrating that the potential softening criterion is not a numerically observed continuation boundary.
At the same time, the map does not directly measure the onset of physical instability.

The plotted quantity $T_{\rm conv}$ depends on nonlinear solver convergence, the warm-start procedure, and the numerical configuration.
A run that reaches $T_{\rm max}=601$ may remain bounded over the available interval, but the map alone cannot establish a physical lifetime of 601, a globally bounded solution, or an oscillon-mediated metastable island.

The most direct conclusion from Fig.~\ref{fig:phi6_continuation_map} is that the lifted interaction produces a broad long-converged region followed by substantially shorter accessible terminal times at selected larger amplitudes.

\subsection{Localization diagnostics}
\label{ssec:phi6_localization}
To examine whether the long-converged solutions exhibit persistent spatial localization, we compare two nearby initial amplitudes at fixed width $\sigma=0.11L$.

The selected cases are
\begin{equation}
A/A_c=1.4
\qquad\text{and}\qquad
A/A_c=1.5.
\end{equation}

We analyze the second-moment radius $\sigma_x(t)$ defined in Eq.~\eqref{eq:sigmax}, together with the field-amplitude histories.

\begin{figure}[H]
\centering
\includegraphics[width=\linewidth]{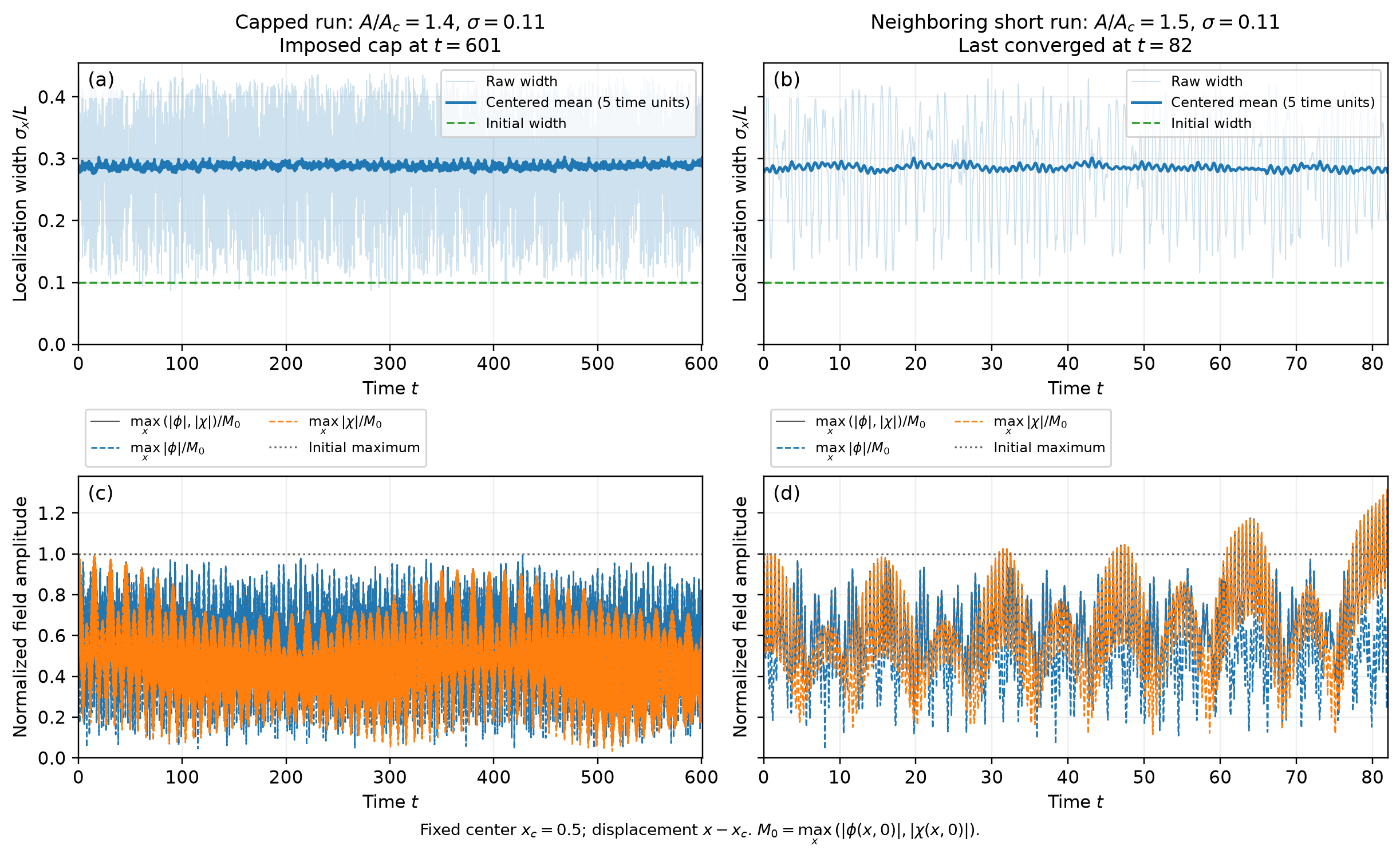}

\caption{Localization radius and normalized field-amplitude histories for two lifted-model runs with $\sigma=0.11L$ and $A/A_c=1.4,1.5$. The smoothed localization radii in both cases are close to $L/\sqrt{12}$, the reference value for uniform spatial weight over a box of length $L$. The selected histories do not show an evidently persistent compact oscillon core. The shorter-continuation case develops larger late-time amplitude excursions.}
\label{fig:phi6_localization_comparison}
\end{figure}

For the initial profile
\begin{equation}
f(x)=
\sech\left(
\frac{x-x_0}{\sigma}
\right),
\end{equation}
the second-moment localization radius on a sufficiently large domain is approximately
\begin{equation}
\sigma_{x,0}
\simeq
\frac{\pi\sigma}{\sqrt{12}}.
\label{eq:sech_initial_radius}
\end{equation}

For $\sigma=0.11L$, this gives
\begin{equation}
\sigma_{x,0}
\simeq0.10L.
\end{equation}

By comparison, a uniform spatial weight on a domain of length $L$ has
\begin{equation}
\sigma_{\rm uniform}
=
\frac{L}{\sqrt{12}}
\simeq0.289L.
\end{equation}

The smoothed radii in the selected evolutions are approximately
\begin{equation}
\sigma_x\simeq
0.28\text{--}0.29L.
\end{equation}
These values are substantially closer to the uniform-weight reference than to the estimated radius of the initial localized profile.

Consequently, Fig.~\ref{fig:phi6_localization_comparison} does not provide positive evidence that the selected long-converged configuration remains a compact, oscillon-like lump.
The result is more consistent with spatially extended field weight over the plotted interval, although a second-moment radius alone cannot establish the detailed spatial structure of the solution.

The long-converged case maintains comparatively bounded field amplitudes over the displayed interval, while the neighboring case exhibits larger late-time excursions.
The distinction between their continuation times therefore cannot be attributed to persistent compact localization on the basis of these diagnostics.
In particular, we do not interpret the broad region of Fig.~\ref{fig:phi6_continuation_map} as a demonstrated oscillon island.

\subsection{Comparison with the quartic model}
\label{ssec:phi6_comparison}
We finally compare the last converged terminal times for the lifted potential and the quartic baseline using the same localized initial-data family.

The comparison uses
$\sigma=0.11L$
and the five
normalized amplitudes
\begin{equation}
A/A_c\in
\{1.3,1.4,1.5,1.6,1.75\}.
\end{equation}

\begin{figure*}[!t]
\centering
\includegraphics[width=0.85\textwidth]{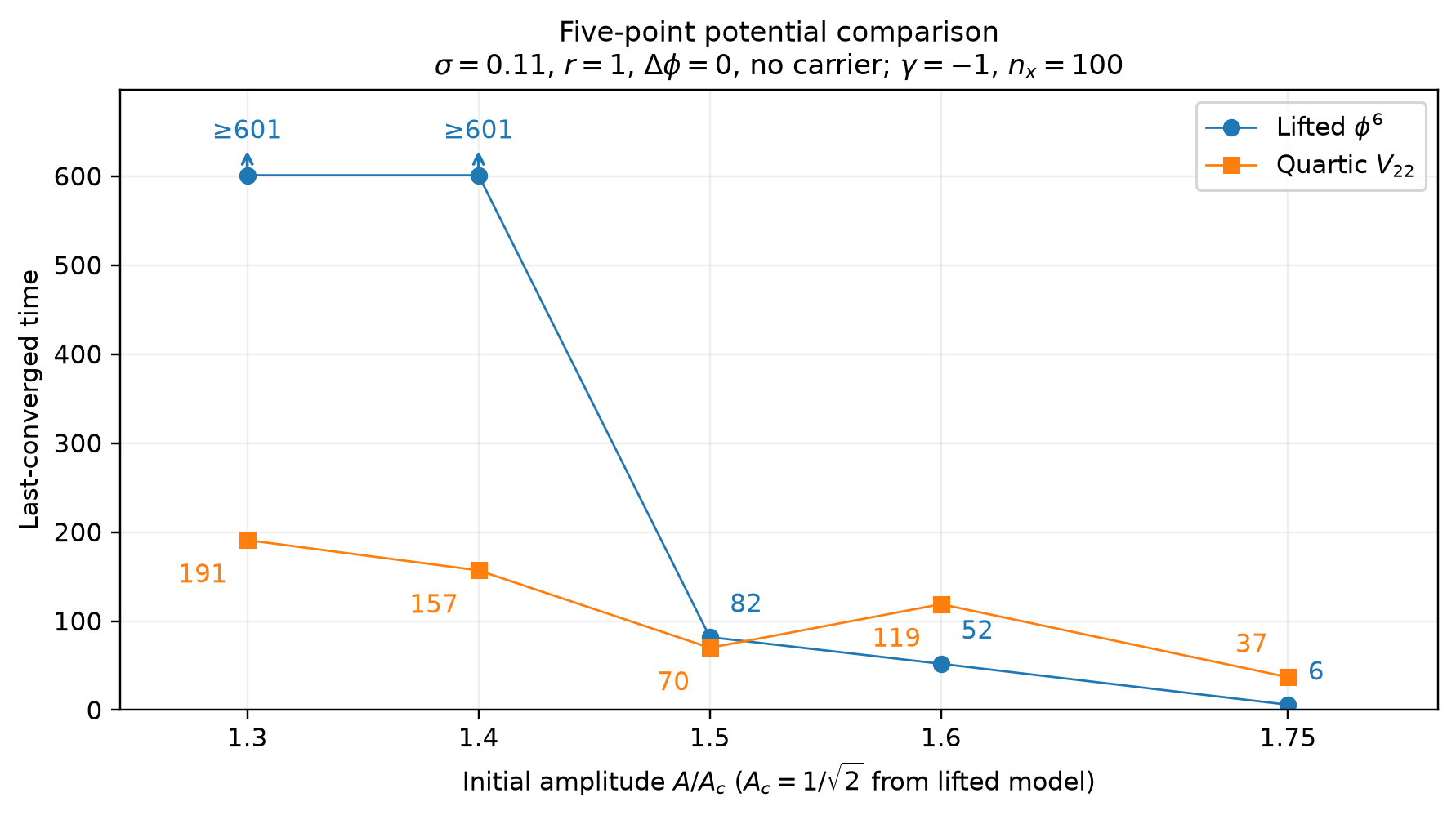}
\caption{Comparison of the last converged terminal time
$T_{\rm conv}$ for the lifted $\phi^6$ and quartic interaction models using the same localized initial-data family.
The lifted model reaches the prescribed cap at the two lowest displayed amplitudes, while its accessible terminal time shortens at larger amplitudes. The comparison involves the complete field equations, including different linearized mass spectra, and therefore cannot isolate stabilization by the sextic term.}
\label{fig:phi6_quartic_comparison}
\end{figure*}

Figure~\ref{fig:phi6_quartic_comparison} shows that the lifted model reaches the computational cap for the two lowest sampled amplitudes, $A/A_c=1.3$ and $A/A_c=1.4$, whereas the corresponding quartic-model calculations terminate earlier.

At larger amplitudes, the lifted-model continuation times shorten and exhibit non-monotonic variation. Neither potential is uniformly dominant at every amplitude.
The comparison thus demonstrates a substantial dependence of the accessible finite-time dynamics on the specified field-theory
model.

However, two limitations prevent interpreting this difference as an isolated physical stabilization effect.
First, the plotted quantity is $T_{\rm conv}$ rather than an independently measured physical instability time.
Second, as established in Sec.~\ref{ssec:phi6_model}, the lifted model modifies the linearized masses as well as the nonlinear potential.
The difference between the curves therefore reflects the combined effect of these changes and their influence on the numerical continuation problem.

The present results support the conclusion that the full interaction model strongly reorganizes the continuation landscape.
They do not establish that oscillon formation is responsible for the difference or that adding a sextic interaction at fixed masses necessarily increases physical stability.
The role of genuinely localized oscillons in ghost-normal systems, and the influence of isolated higher-order interactions at matched linear dispersion relations, remain questions for future studies.

\section{Conclusion}
\label{sec:conclusion}

We have examined the finite-time classical dynamics of nonlinear scalar fields with one normal and one ghost kinetic sector.
The simulations reveal a structured dependence of the observed behavior on initial amplitude, spectral distribution, waveform localization, and interaction model. They do not support a universal immediate classical runaway, but neither do they establish global stability.

The relevant numerical distinction is between strong growth observed within a resolved evolution, bounded motion over the sampled interval, and termination caused by the nonlinear algebraic solver.
For the quartic interaction $V_{22}=\lambda_{22}\phi^2\chi^2$, the amplitude--wavenumber and Gaussian-width maps show a consistent qualitative ordering.
Larger amplitudes tend to cross growth thresholds earlier, while higher characteristic wavenumbers or narrower localized packets tend to remain below those thresholds for longer times.
Random-phase colored spectra sharpen this result by varying the spectral tilt at fixed field RMS amplitude.
The UV-biased cases generally remain uncrossed for longer than the
IR-biased cases, but the measured onset times also display substantial realization-to-realization variation.

These comparisons establish spectral organization of finite-time threshold crossing, not a universal instability timescale or a control experiment at fixed total Hamiltonian.
The additional initial-data families delimit the interpretation.
The phase-correlated maps suggest stronger dependence on relative waveform alignment in co-propagating configurations, although solver-limited cases and unequal threshold choices preclude a sharp quantitative comparison with the counter-propagating survey.

Initially compact and time-symmetric perturbations exhibit a broad region without observed threshold crossing.
Because the simulations use periodic spatial boundaries, this observation is a finite-time analogue of compact-data questions, not numerical evidence for scattering to infinity.

Several nonlinear controls also show that loss of solver convergence can occur while the field amplitudes remain bounded, underscoring the importance of separating algorithmic termination
from physical growth.

Time-resolved Gaussian diagnostics provide a more direct connection between spectral behavior and field growth. In one rapidly growing case, low-wavenumber Fourier power increases before the large field-amplitude excursion, whereas a neighboring long-lived case does not exhibit a comparable increase.
This is evidence for a spectral precursor in the selected evolution. Its causal role and generality remain open: the ratio of field Fourier powers is not a signed Hamiltonian-energy flux, and ordinary mode coupling alone is not sufficient to identify a cascade.

The lifted $\phi^6$ model gives an additional demonstration of the sensitivity of long-time numerical evolution to the interaction model. The amplitude--width scan contains a broad region whose nonlinear solve reaches the imposed terminal-time cap, followed by much shorter accessible continuation times at larger amplitudes.
However, the lifted quadratic term also changes the linearized mass spectrum, so this comparison cannot isolate the effect of sextic self-interaction.
Furthermore, the localization radii observed in two selected cases are not indicative of persistent compact oscillon cores.
The supported conclusion is therefore potential-dependent continuation behavior, rather than demonstrated oscillon-mediated stabilization.
More generally, the study shows how a metastability landscape can be characterized without conflating operational threshold times, censored runs, and solver failures.

Future work could determine whether the observed spectral precursors correspond to a mode-resolved transfer of signed energy, extend the analysis to scattering on unbounded domains, and isolate the influence of higher-order interactions at matched linear mass spectra.
These questions arise naturally from the current results, but do not weaken their main conclusion:
In classical ghost-normal wave systems, finite-time dynamics is highly structured and cannot be inferred from the sign of the kinetic term alone.
Mapping this metastability landscape provides a practical way to understand when ghost-normal systems undergo prompt runaway, when they remain long lived, and how nonlinear structure can delay---though not necessarily eliminate---ghost-driven instability.

\paragraph*{Acknowledgements.}
We thank Aaron Held for helpful discussions and directions. 
This work used resources provided by SDSU's Innovator HPC Cluster, MSU's High-Performance Computing Cluster, and  the LANL Darwin testbed. Darwin is a research testbed/heterogeneous cluster funded by the Computational Systems and Software Environments subprogram of ASC program. LANL is operated by Triad National Security, LLC, for the National Nuclear Security Administration of the U.S.DOE  (Contract No. 89233218CNA000001). This work is authorized for unlimited release under LA-UR-25-31553.

\paragraph*{Artificial intelligence assistance.} 
The authors acknowl edge the limited use of artificial intelligence tools, including ChatGPT (OpenAI), during the preparation of this manuscript. These tools were used for spelling and grammar checks, LaTeX editing and debugging, and improving the clarity and readability of some text. They were also used to assist with code debugging and code
development, mostly in generating scripts to analyze and plot the data. All main methods/modeling developments and simulations are performed by the authors. The authors reviewed and verified all AI-assisted material and take full responsibility for the scientific content, calculations, code, interpretations, and conclusions presented in this work

\bibliographystyle{unsrt}
\bibliography{refs}

\appendix
\section{Space-Time Finite Element Method}
\label{sec:appx:sFEM}

A space-time FEM uses continuous approximation functions in both space and time. We follow the discretization scheme in~\cite{Anderson2007, French1996}. 
In this scheme, space and time are discretized together for the entire domain using a finite element space which does not discriminate between space and time basis functions. 
In this way, we solve the entire spacetime slab simultaneously.
This is different from other methods that discretize space and then step through time slabs, \cite{Deffayet2025, Cherubini2005(1), Cherubini2005(2), Cao2018}.
The space-time FEM discretization scheme has been used to investigate interesting phenomena in general relativity, \cite{Lim2020}. The method has also been utilized in engineering applications \cite{Kim2001, Guddati1999, Dyniewicz2012} and to solve the wave equation to display some interesting numerical finding, \cite{Anderson2007, Zank2025}.

Thus, the work done here also serves to motivate the use of the space-time FEM discretization on nonlinear problems in physics applications.
Our setup employs a structured mesh using PETSc's DMDA framework, \cite{petsc-web-page}, which provides a grid with uniform spacing between nodes. However, we note that it is sometimes desirable to employ an unstructured mesh on the domain where the spacing between nodes varies according to the application, \cite{Schwing2013, Xu2023, Gou2016}. (This option is also available within PETSc by employing their DMPlex framework.) 

For example, consider an engineer using the FEM to model a bridge. If the goal of the work is to better understand the stress distribution near the center of the bridge, then an unstructured mesh may be employed. In this way, discretized domain elements near the center will be smaller than those at the edges. This will allow for sharper results about the area of interest while still allowing the engineer to model contributions from the entire system.

By analogy, adaptive or unstructured meshes are sometimes useful in modeling phenomena occurring at or near massive objects that greatly warp space-time, \cite{Bryan2014, Xu1997, Evans2005}. This would allow for the simulation to better capture the dynamics of the system near the high curvature regions while still taking into account any other phenomena in the model. While this is beyond the scope of the present work, it motivates the use of space-time FE discretizations in exploring other interesting physical phenomena.

\subsection{Setting Up The Weak Form}
Since the original system is second order in time, it is necessary for the weak formulation to introduce auxiliary variables such that our system can be reduced to one that is first order in time. 

Let $u = \partial_t \phi$ and $v = \partial_t\chi$. From Eqns.~\ref{eq:phi_strong_EOM}~and~\ref{eq:chi_strong_EOM}, we have
\begin{align}
\partial_t u - \nabla^2\phi + m_{\phi}^2\phi + \partial_\phi V &= 0, \\
\frac{\partial \phi}{\partial t} - u &= 0, \\
\partial_t v - \nabla^2\chi + m_{\chi}^2\chi + \gamma\partial_\chi V &= 0, \\
\frac{\partial \chi}{\partial t} - v &= 0.
\label{eqn:eom:1st_order_storng}
\end{align}

Then, we multiply all equations by a test function (rectangular basis) and integrate over the spacetime domain. Integration by parts in space allows us to move a derivative to the test function and utilize our spatially periodic BC to remove the spatial boundary term integrals. Thus, with spatial and temporal domains of $\Omega$ and T, the weak form of the EOM is given by:
\begin{align}
    K_1 &= \int_{\Omega,T}\!\Bigl(u_t\Psi+\nabla\phi\nabla\Psi
                     \nonumber\\
        &\qquad +m_{\phi}^2\phi\Psi+\partial_\phi V\,\Psi\Bigr)
          \,d\Omega\,dt,\nonumber\\
    K_2 &= \int_{\Omega,T}(\phi_t\Psi-u\Psi)\,d\Omega\,dt,
           \nonumber\\
    G_1 &= \int_{\Omega,T}\!\Bigl(v_t\Psi+\nabla\chi\nabla\Psi
                     \nonumber\\
        &\qquad +m_{\chi}^2\chi\Psi+\gamma\partial_\chi V\,\Psi\Bigr)
          \,d\Omega\,dt,\nonumber\\
    G_2 &= \int_{\Omega,T}(\chi_t\Psi-v\Psi)\,d\Omega\,dt.
    \label{eq:Weak_EOM}
\end{align}

There are several different ways to define basis functions.
Here, the basis function is described by rectangular elements. With $\xi = x/h_x$ and $\tau = t/h_t$, where $h_x$ and $h_t$ refer to the distance between nodes in the mesh, we define the rectangular basis functions for the (1 + 1) case:
\begin{align*}
    \Psi_1(\xi, \tau) &= (1 - \xi)(1 - \tau), \\
    \Psi_2(\xi, \tau) &= \xi(1 - \tau), \\
    \Psi_3(\xi, \tau) &= \xi \tau, \\
    \Psi_4(\xi, \tau) &= (1 - \xi) \tau.
\end{align*}

For the (2 + 1) scenario, we follow the same formulation.
With $\xi = x/h_x$, $\tau = t/h_t$, and $\zeta = y/h_y$, we define
the trilinear basis functions:
\begin{align*}
\Psi_1(\xi,\tau,\zeta) &= (1 - \xi)(1 - \tau)(1 - \zeta), \\
\Psi_2(\xi,\tau,\zeta) &= \xi(1 - \tau)(1 - \zeta), \\
\Psi_3(\xi,\tau,\zeta) &= \xi \tau (1 - \zeta), \\
\Psi_4(\xi,\tau,\zeta) &= (1 - \xi)\tau (1 - \zeta), \\
\Psi_5(\xi,\tau,\zeta) &= (1 - \xi)(1 - \tau)\zeta, \\
\Psi_6(\xi,\tau,\zeta) &= \xi(1 - \tau)\zeta, \\
\Psi_7(\xi,\tau,\zeta) &= \xi\tau\zeta, \\
\Psi_8(\xi,\tau,\zeta) &= (1 - \xi)\tau\zeta.
\end{align*}

Using these basis functions, the element stiffness matrix is assembled.  

\subsection{Nonlinear Solver: PETSc SNES}
\label{sec: snes explanation}
In order to computationally handle the nonlinearity of the PDE system, we utilize PETSc's SNES library. The library contains methods, like Newton's method with line search, for solving nonlinear equations of the form
\begin{equation*}
    \mathcal{F}(U) = 0,
\end{equation*}
where $\mathcal{F}$ is the nonlinear differential operator, and $U$ is the solution vector.

The two main user-created components of the code are the $\texttt{FormResidual()}$ and \texttt{FormJacobian()} functions.
$\texttt{FormResidual()}$ takes as input an approximate guess of the solution vector $U$. 
Here, the solution is inserted into the weak form of the EOM. 
Over each element, the integration is done with the pre-calculated element stiffness matrices for the linear terms, and Gaussian quadrature for the nonlinear terms.

This process creates a residual vector, $\mathcal{F}(U) = R$.
If the norm of the residual vector is near zero, then we have solved $\mathcal{F}(U) \approx 0$, and we are done.
If the current $U$ is not acceptable, then we aim to update $U$ in a way that minimizes $\mathcal{F}(U)$.

To do this, \texttt{FormJacobian()} takes in the current iterate $U_k$ and then implements and returns $J =\frac{\partial \mathcal{F}}{\partial U}|_{U_k}$. Then, a perturbation of $U_k$ that should decrease the residual $R_k$ is found by solving $J \delta U = - R_k$. Now, we find a new guess by replacing $U_k$ with $U_k + \delta U_k$ and repeat the process.

The only parts of this process that we manually implement are the creation of the residual vector, $R$, and the Jacobian matrix, $J$, that are made by solving the weak EOM over each element. 
These routines will necessarily change from user to user depending on their choice of numerical discretization.
PETSc's SNES call handles all the other operations (checking the norm of the residual, solving $J \delta U_k = -R_k$, and updating the iterate according to a line search).
See table \ref{tab: snes process} for a step-by-step process.
\begin{table*}[!t]
\[
\begin{array}{l}
\hline
\textbf{SNES Algorithm (Newton with line search)} \\ \hline
\textbf{Input:} \text{ initial guess } U_0 \text{ (seeded from initial conditions)} \\
\textbf{Repeat for } k=0,1,2,\dots \\
\quad 1.\; R_k \gets F(U_k) \quad \texttt{(FormResidual())} \\
\quad 2.\; \text{If } \|R_k\| \le \text{rtol},\; \text{stop (converged)} \\
\quad 3.\; J_k \gets \dfrac{\partial F}{\partial U}(U_k) 
\quad \texttt{(FormJacobian())} \\
\quad 4.\; \text{Solve } J_k\,\delta U_k = -R_k \quad \text{(KSP/PC linear solve, PETSc)} \\
\quad 5.\; \text{Choose step } \alpha \in (0,1] \text{ (line search, PETSc)} \\
\quad 6.\; U_{k+1} \gets U_k + \alpha\,\delta U_k \quad \text{(PETSc)}\\
\hline
\end{array}
\]
\caption{}
\label{tab: snes process}
\end{table*}

\label{sec:appx:stiff_mat}
An important step in any FEM is determining the element stiffness matrices.
These matrices are found by solving the weak form of the equations of motion, Eq. \ref{eq:Weak_EOM}.
Solving the space-time FEM for the (1 + 1) case results in 4x4 matrices, the (2 + 1) case finds 8x8 matrices, and the 3+1 case has 16x16 matrices. 

The weak equations are displayed in equation \ref{eq:Weak_EOM}, but in equations \ref{eq: K1}, \ref{eq: K2}, \ref{eq: G1}, and \ref{eq: G2}  we display them in their (1 + 1) form.
\begin{align}
    K_1 &= \int_{X,T}\!\Bigl(u_t\Psi+\phi_x\Psi_x
                   \nonumber\\
         &\qquad+m_\phi^2\phi\Psi+V_\phi\Psi\Bigr)\,dx\,dt
                   \label{eq: K1}\\
    K_2 &= \int_{X,T}(\phi_t\Psi-u\Psi)\,dx\,dt
                   \label{eq: K2}\\
    G_1 &= \int_{X,T}\!\Bigl(v_t\Psi+\chi_x\Psi_x
                   \nonumber\\
         &\qquad+m_\chi^2\chi\Psi+\gamma V_\chi\Psi\Bigr)\,dx\,dt
                   \label{eq: G1}\\
    G_2 &= \int_{X,T}(\chi_t\Psi-v\Psi)\,dx\,dt
                   \label{eq: G2}
\end{align}
Recall that $\Psi$ in the above equations is the test function.
The test function is a piecewise polynomial function and can take on many different forms.
For our purposes, a rectangular basis function is used.
Once the rectangular element basis function is determined, it will be inserted into the weak equations.
These integrals will be solved over each element and element stiffness matrices will be produced.

\subsection{Element Basis Function}
The rectangular element basis function is a piecewise polynomial with four terms,
\begin{equation}
    \Psi_i(x,t) = b_{1,i}x + b_{2,i}t + b_{3,i}xt + b_{4,i}.
    \label{eq: phi}
\end{equation}
To solve for the unknown coefficients in Eq. \ref{eq: phi}, we consider the four corners of a rectangle in the $(x,t)$ plane, see Fig. \ref{fig: basisElementFunc}.
\begin{figure}[H]
\centering
    \begin{tikzpicture}[scale = 2]
    \draw[step=1cm,gray] (0,0) grid (1,1);
    \node at (-0.5, 0.0) {(0,0)};
    \fill (0,0) circle (3pt);
    \node at (1.5, 0.0) {($h_x$,0)};
    \fill (1,0) circle (3pt);
    \node at (-0.5, 1) {(0,$h_t$)};
    \fill (0,1) circle (3pt);
    \node at (1.5, 1) {($h_x$,$h_t$)};
    \fill (1,1) circle (3pt);
    \end{tikzpicture}
\caption{}
\label{fig: basisElementFunc}
\end{figure}
In figure \ref{fig: basisElementFunc}, $h_x$ and $h_t$ represent the step size.
They can be thought of as the distance between spatial and temporal nodes. They are defined by
\begin{align*}
    h_x &= \frac{x_{\text{final}} - x_{\text{initial}}}{n_x} \\
    h_t &= \frac{t_{\text{final}} - t_{\text{initial}}}{n_t - 1} \\
\end{align*}
where $n_x$ and $n_t$ are the number of spatial and temporal nodes in the entire domian.
Thus, $n_x$ and $n_t$ control the resolution of the simulation; as the number of nodes grows, the distance between nodes shrinks, and the domain becomes more refined.
It should be noted that the ($n_x$) and ($n_t -1$) terms in the denominators of the above equations is not a mistake.
It must be this way since we implement periodic spatial boundaries but impose no such condition temporally.

Since we are currently looking for the element stiffness matrix, only  one `element' is examined in the entire domain.
The element can be visualized as has been done in Fig. \ref{fig: basisElementFunc} with each corner representing a node designated by the subscript $i$.
The coordinates $(0, 0), \ (h_x,0), \ (0,h_t), \ (h_x,h_t)$ will be inserted into Eq. \ref{eq: phi} so that the $b$ coefficients can be found. This results in four separate $\Psi$ equations. To determine the basis function, we require that
\begin{equation}
    \Psi_{i} = 
    \begin{cases} 
        1, & \text{if } j = i \\
        0, & \text{if } j \neq i 
    \end{cases}
    \nonumber
\end{equation}
where $j$ is represented as is described in Table \ref{tab: table}.
\begin{table}
    \centering
    \caption{Representation of $j$.}
    \begin{tabular}{c c c}
    \hline
    $j$ & $x$ & $t$ \\
    \hline
    1 & 0 & 0 \\
    \hline
    2 & 0 & $h_t$ \\
    \hline
    3 & $h_x$ & 0 \\
    \hline
    4 & $h_x$ & $h_t$ \\
    \hline
    \end{tabular}
    \label{tab: table}
\end{table}

This results in four polynomials, the element basis functions
\begin{align}
    \Psi_1 &= -\frac{1}{h_x}x -\frac{1}{h_t}t + \frac{1}{h_x h_t}xt + 1\text{,} \nonumber \\
    \Psi_2 &= \frac{1}{h_x}x - \frac{1}{h_x h_t}xt, \nonumber \\
    \Psi_3 &= \frac{1}{h_x h_t}xt,  \nonumber \\
    \Psi_4 &= \frac{1}{h_t}t - \frac{1}{h_x h_t}xt. \nonumber
\end{align}

\subsection{Element Stiffness Matrices}
Now that we have the basis functions, the element stiffness matrices can be determined.
We will explain the set up of the element stiffness matrices, display the matrices, and then discuss practical implementation into the simulation.

In the two dimensional (1 + 1) case, the weak form of the equations of motion, equations \ref{eq: K1}, \ref{eq: K2}, \ref{eq: G1}, and \ref{eq: G2}, are double integrals over one space-time element in the domain,
\begin{align*}
    K_1 &= \int_0^{h_x}\!\int_0^{h_t}\!\Bigl(u_t\Psi+\phi_x\Psi_x
                   \\
         &\qquad +m_\phi^2\phi\Psi+V_\phi\Psi\Bigr)\,dx\,dt\\
    K_2 &= \int_0^{h_x}\!\int_0^{h_t}(\phi_t\Psi-u\Psi)\,dx\,dt\\
    G_1 &= \int_0^{h_x}\!\int_0^{h_t}\!\Bigl(v_t\Psi+\chi_x\Psi_x
                   \\
         &\qquad +m_\chi^2\chi\Psi+\gamma V_\chi\Psi\Bigr)\,dx\,dt\\
    G_2 &= \int_0^{h_x}\!\int_0^{h_t}(\chi_t\Psi-v\Psi)\,dx\,dt 
\end{align*}

For demonstration purposes, we will put together the element stiffness matrices pertaining to the three linear terms in the weak EOM.
The linear terms are the ones that do not include contributions from $V$ or its partial derivatives. 
This will result in three $(4 \times 4)$ element stiffness matrices each of the form
\begin{equation}
A = 
\begin{bmatrix}
a_{1,1} & a_{1,2} & a_{1,3} & a_{1,4} \\
a_{2,1} & a_{2,2} & a_{2,3} & a_{2,4} \\
a_{3,1} & a_{3,2} & a_{3,3} & a_{3,4} \\
a_{4,1} & a_{4,2} & a_{4,3} & a_{4,4}
\end{bmatrix}
.
\end{equation}

The Time, Space, and Standard element stiffness matrices come from solving
\begin{align*}
    \text{Time}_{i,j} &= \int_0^{h_t}\!\int_0^{h_x} (\partial_t\Psi_i)\Psi_j\,dx\,dt \\
    \text{Space}_{i,j} &= \int_0^{h_t}\!\int_0^{h_x} (\partial_x\Psi_i)(\partial_x\Psi_j)\,dx\,dt \\
    \text{Standard}_{i,j} &= \int_0^{h_t}\!\int_0^{h_x} \Psi_i\Psi_j\,dx\,dt.
\end{align*}
Since these terms are linear in the weak EOM, we are able to calculate their element stiffness matrices one time.
This can be done by hand, but we choose to utilize the SymPy package in Python that is capable of handling symbolic integration.
\[
\begin{aligned}
\text{Time} &=
\dfrac{h_x}{12}
\begin{bmatrix}
 -2 & -1 & 1 & 2 \\
 -1 & -2 & 2 & 1 \\
 -1 & -2 & 2 & 1 \\
 -2 & -1 & 1 & 2
\end{bmatrix},
\\[8pt]
\text{Space} &= 
\dfrac{h_t}{6h_x}
\begin{bmatrix}
 2 & -2 & -1 & 1 \\
 -2 & 2 & 1 & -1 \\
 -1 & 1 & 2 & -2 \\
 1 & -1 & -2 & 2
\end{bmatrix},
\\[8pt]
\text{Mass} &=
\dfrac{h_t h_x}{36}
\begin{bmatrix}
 4 & 2 & 1 & 2 \\
 2 & 4 & 2 & 1 \\
 1 & 2 & 4 & 2 \\
 2 & 1 & 2 & 4
\end{bmatrix}.
\end{aligned}
\]
These matrices are hard-coded into the simulation and saved once, since they never change.

The nonlinear terms of the weak EOM must be dealt with a little differently, but the gist is the same.
Instead of being able to find the element stiffness matrices once at the beginning of the simulation, we must use a Gaussian quadrature routine to approximate the integrals of the potential terms over each element.
This is done inside the \texttt{FormResidual()} and \texttt{FormJacobian()} routines.

\section{Dynamical signatures of runaway and solver-limited evolution}
\label{sec:appx:runaway-diagnostics}

An important numerical distinction in this work is that loss of convergence of the nonlinear solver does not, by itself, imply a physical ghost instability.
We therefore classify representative evolutions using both field amplitudes and the Hamiltonian decomposition introduced in Sec.~\ref{sec:model}.

For a genuine ghost-driven runaway, the characteristic behavior is not necessarily growth of the conserved total Hamiltonian.  Rather, because the normal and ghost sectors carry energy with opposite signs, nonlinear energy exchange can generate increasingly large positive- and negative-energy sector contributions while preserving the total energy,
\begin{equation}
    H(t)
    =
    H_\phi(t)+H_\chi(t)+H_{\rm int}(t)
    \simeq H(0).
\end{equation}
A representative runaway-like configuration therefore exhibits simultaneous growth of the field amplitude and increasingly large compensating excursions
of $H_\phi$ and $H_\chi$.

This behavior should be contrasted with a solver-limited evolution.  In such a case, SNES loses convergence at a finite terminal time even though the field amplitudes remain bounded and no corresponding growth of the sector energies is
observed.  We do not interpret such a termination as evidence of physical runaway.

\begin{figure*}[!t]
\begin{center}
    \begin{subfigure}[t]{0.32\textwidth}
        \centering
        \includegraphics[width=\linewidth]
        {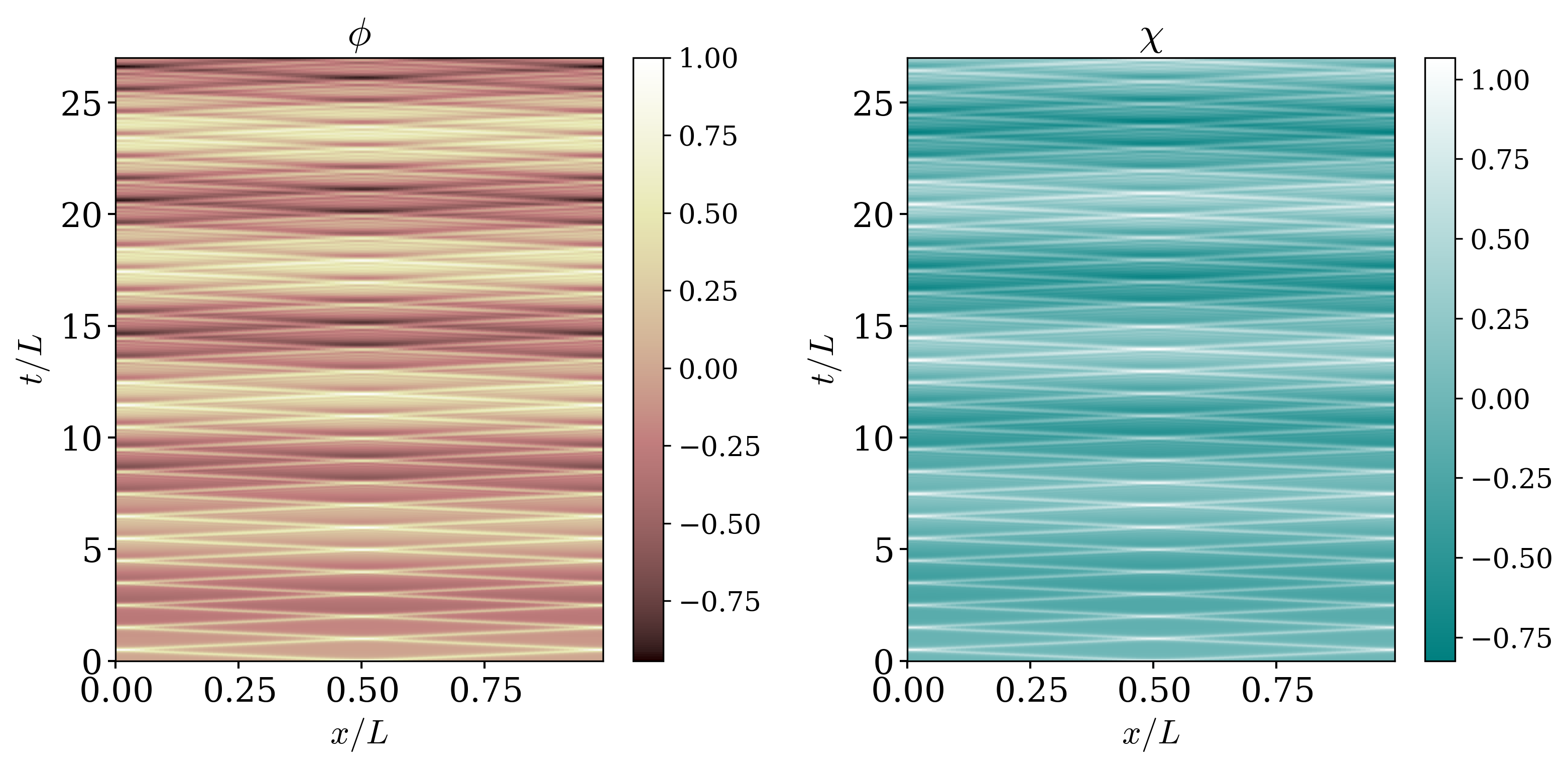}
        \includegraphics[width=\linewidth]
        {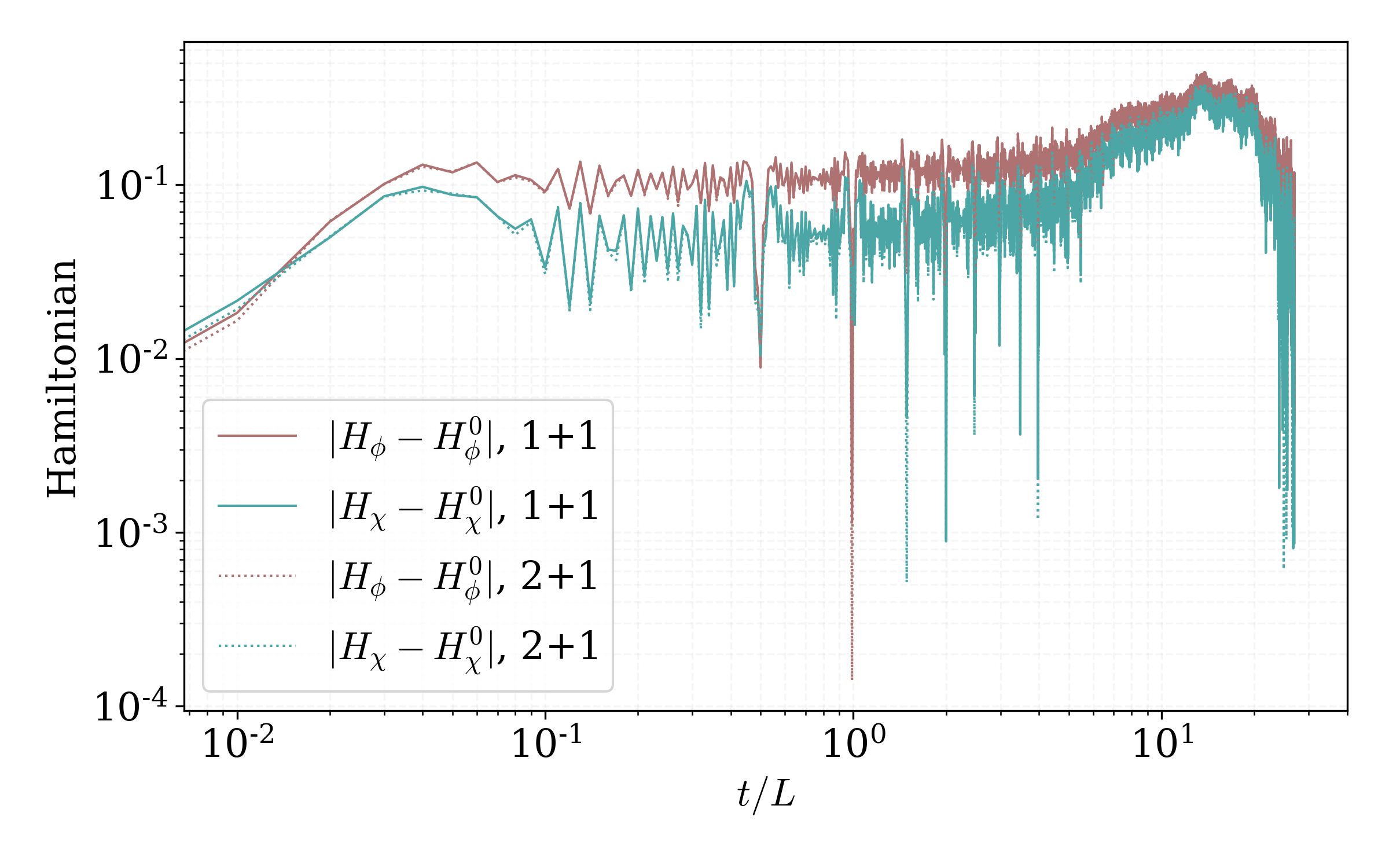}
        \caption{$\gamma=-1$, $\lambda=1$: interacting ghost case.}
        \label{fig:diag-runaway}
    \end{subfigure}
    \hfill
    \begin{subfigure}[t]{0.32\textwidth}
        \centering
        \includegraphics[width=\linewidth]
        {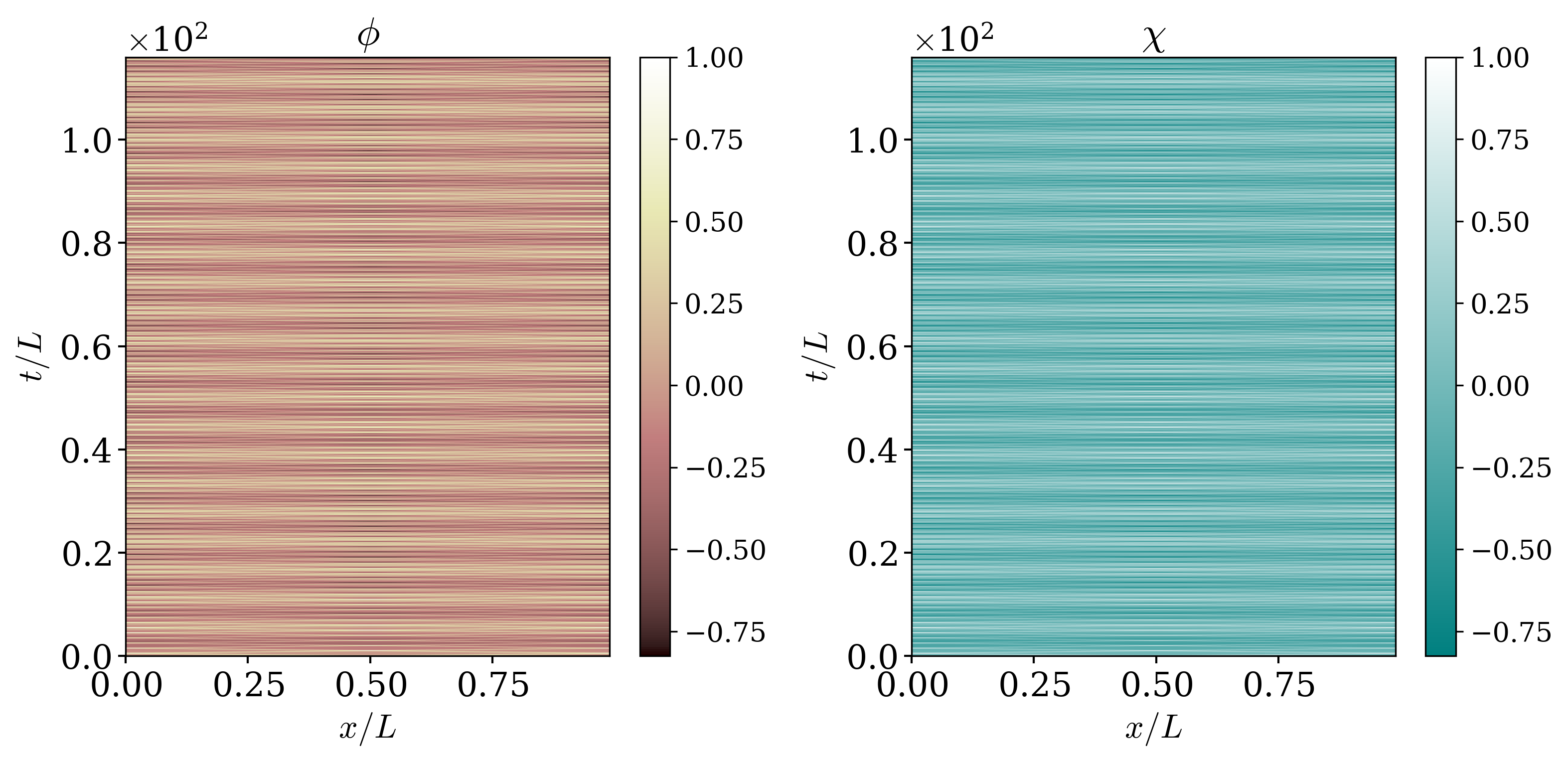}
        \includegraphics[width=\linewidth]
        {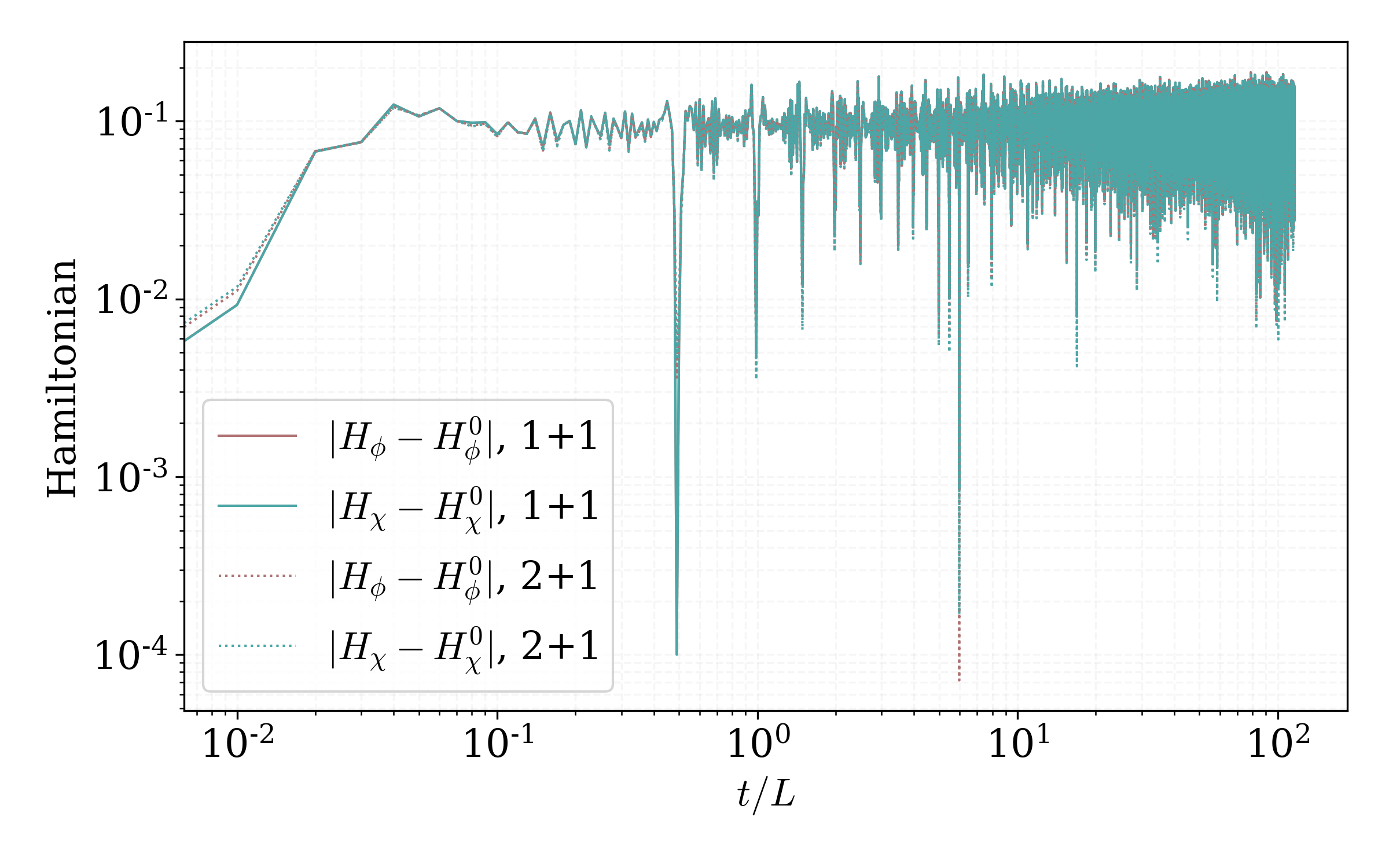}
        \caption{$\gamma=+1$, $\lambda=1$: solver-limited nonlinear case.}
        \label{fig:diag-solver}
    \end{subfigure}
    \hfill
    \begin{subfigure}[t]{0.32\textwidth}
        \centering
        \includegraphics[width=\linewidth]
        {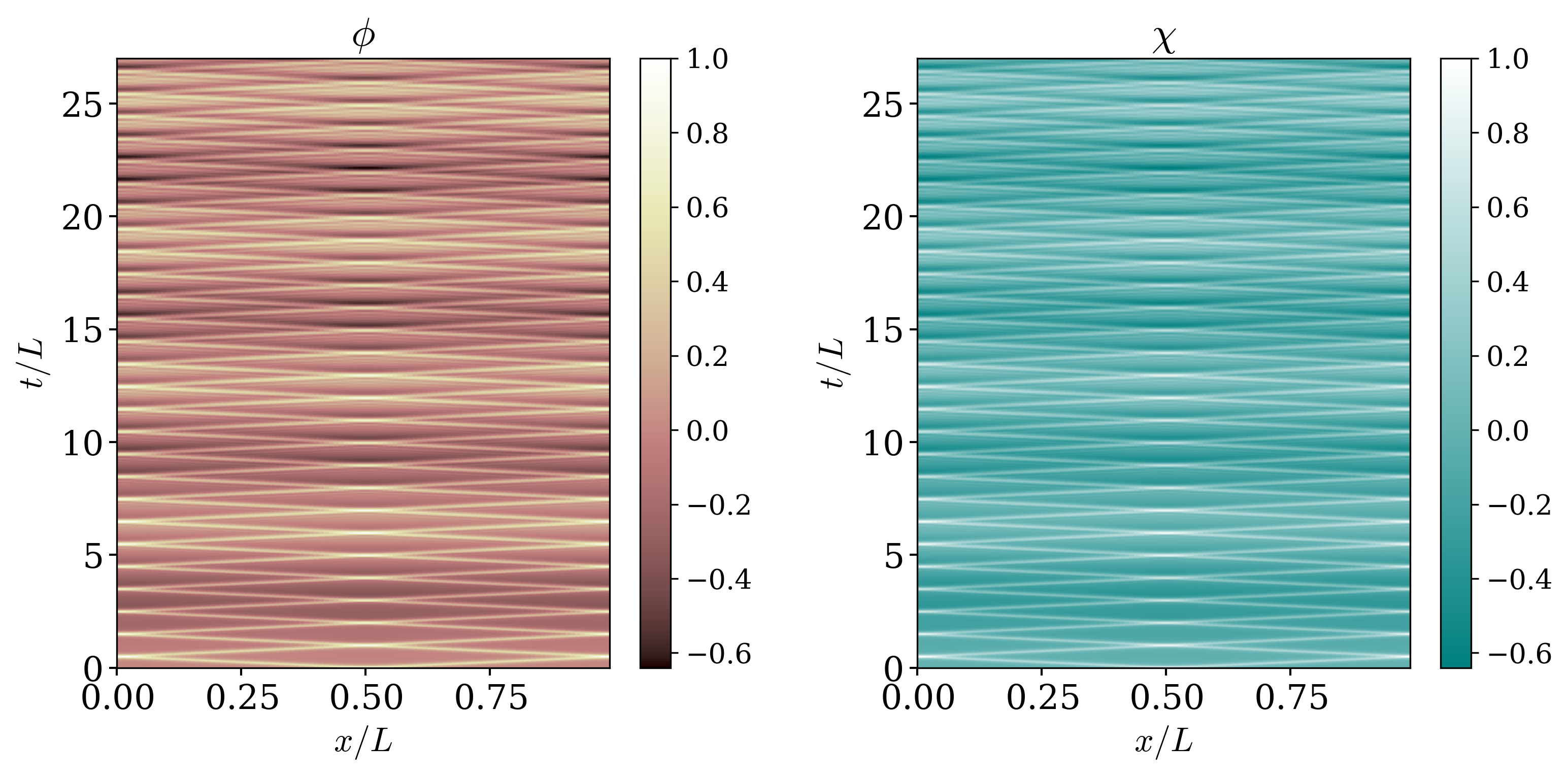}
        \includegraphics[width=\linewidth]
        {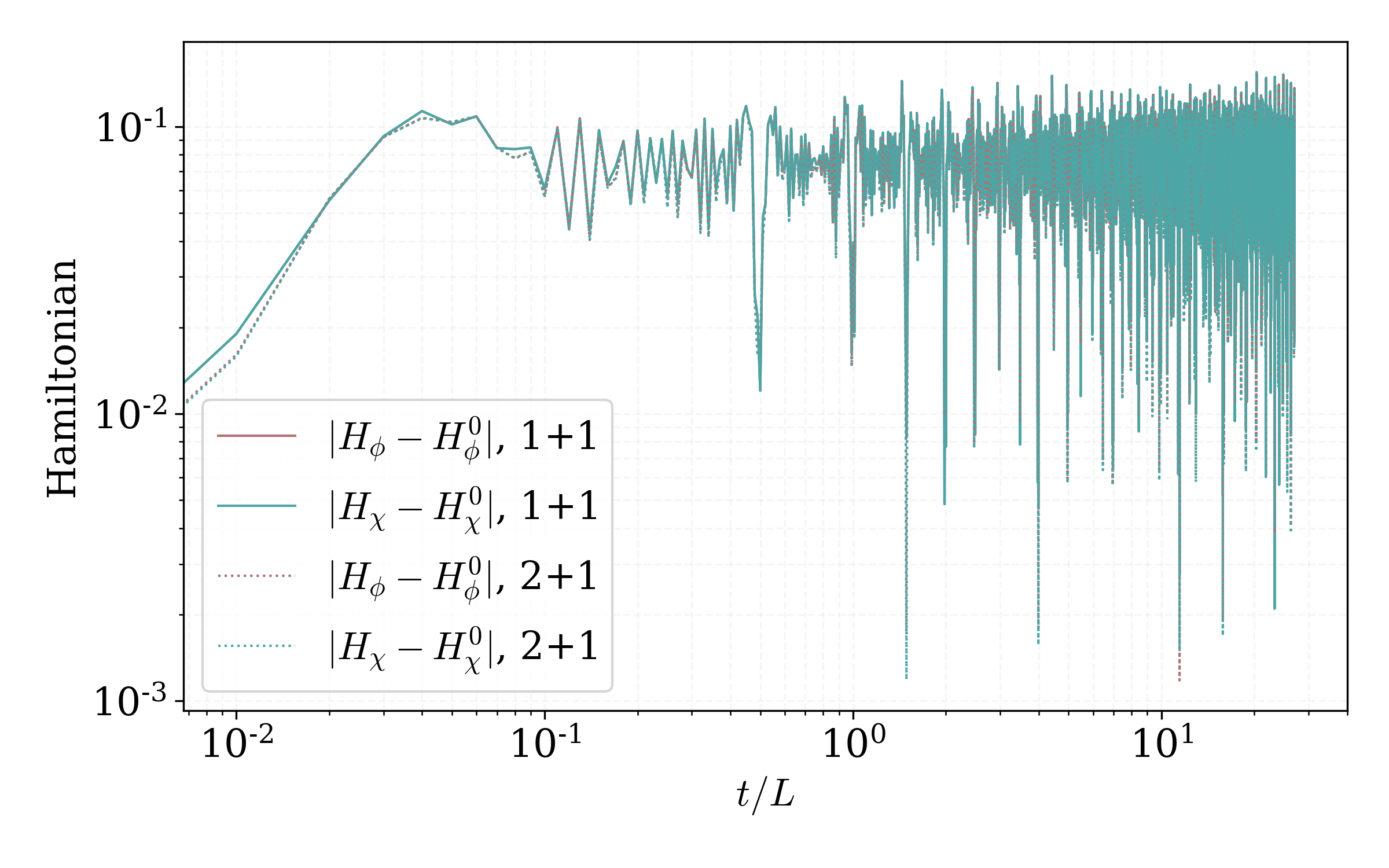}
        \caption{$\gamma=-1$, $\lambda=0$: uncoupled ghost control.}
        \label{fig:diag-control}
    \end{subfigure}
    \end{center}
    \caption{\justifying Representative dynamical signatures used to distinguish physical ghost-driven growth from numerical solver failure.  In the interacting ghost case, field growth is accompanied by increasingly large positive- and negative-energy sector contributions.  In contrast, the solver-limited non-ghost case terminates without comparable amplitude or sector-energy growth, while removing the nonlinear interaction leaves the ghost system bounded over the comparison interval.  The total Hamiltonian is monitored separately to verify that large sector-energy excursions are associated with compensating energy exchange rather than uncontrolled numerical energy drift.}
    \label{fig:runaway-classification}
\end{figure*}

Figure~\ref{fig:runaway-classification} illustrates three representative situations: an interacting ghost system exhibiting runaway-like growth, a non-ghost nonlinear evolution that becomes solver-limited while remaining
bounded, and a ghost system with the nonlinear coupling removed that remains well behaved over the comparison interval.

\section{Convergence and validation}
\label{sec:appx:conv}
\subsection{Manufactured solution tests}\label{ssec: MMS 1+1}
We consider the following coupled partial differential equations (PDEs) in (1 + 1) dimensions:
\begin{align}
  \phi_{tt} - \phi_{xx} + \phi + V_\phi(\phi,\chi) &= 0, \label{eq:phi} \\
  \chi_{tt} - \chi_{xx} + \chi - V_\chi(\phi,\chi) &= 0, \label{eq:chi}
\end{align}
where the highly nonlinear potential terms are defined as
\begin{align}
  V_\phi(\phi,\chi) &= -2\lambda \, \phi\, \frac{A}{B^{1.5}} \nonumber \\ 
  V_\chi(\phi,\chi) &= 2\lambda \, \chi\, \frac{C}{B^{1.5}}  \nonumber
\end{align}
with
\begin{align}
    A &= \phi^2 - \chi^2 + 1 \nonumber \\
    B &= \left(\phi^2 - \chi^2 - 1\right)^2 + 4\phi^2, \nonumber\\
    C &= \phi^2 - \chi^2 - 1.\nonumber
\end{align}
Here, $\lambda$ is a parameter controlling the strength of the nonlinear potential. It is set to $1.0$ for all testing in this section.

\subsubsection{Problem Setup and Goal}

The aim of our numerical implementation is to solve equations \eqref{eq:phi}--\eqref{eq:chi} using a spacetime finite element method (FEM) within PETSc with a nonlinear solver (SNES). The code uses a formulation in which each spacetime node carries four degrees of freedom:
\begin{itemize}
  \item $\phi(x,t)$ and an auxiliary field $u(x,t)$, which is related to the time derivative of $\phi$, and
  \item $\chi(x,t)$ and an auxiliary field $v(x,t)$, which is related to the time derivative of $\chi$.
\end{itemize}

To verify our implementation, we use the \emph{method of manufactured solutions} (MMS). The idea behind MMS is to choose a smooth (and nontrivial) manufactured solution for $\phi(x,t)$ and $\chi(x,t)$, compute all the required derivatives, and then add appropriate forcing functions so that the manufactured solution exactly satisfies the modified PDEs. In our case, we modify the equations to
\begin{align}
  \phi_{tt} - \phi_{xx} + \phi + V_\phi(\phi,\chi) &= F_\phi(x,t), \label{eq:phi_force} \\
  \chi_{tt} - \chi_{xx} + \chi - V_\chi(\phi,\chi) &= F_\chi(x,t). \label{eq:chi_force}
\end{align}
The goal is to compute these forcing functions for $\phi$ and $\chi$ ($F_{\phi}\text{ and }F_{\chi}$) so that when the manufactured solution is substituted into the left-hand side of \eqref{eq:phi_force} and \eqref{eq:chi_force}, the LHS exactly matches $F_\phi(x,t)$ and $F_\chi(x,t)$.

\subsubsection{Manufactured Solution}

We choose the following manufactured solutions:
\begin{align}
  \phi(x,t) &= 1 + \cos\bigl(\pi x\bigr) \cos\bigl(\pi t\bigr), \label{eq:phi_manuf} \\
  \chi(x,t) &= 1 + \cos\bigl(\pi x\bigr) \sin\bigl(\pi t\bigr). \label{eq:chi_manuf}
\end{align}
These functions are chosen because they are smooth, nontrivial, and spatially periodic on symmetric spatial domains.

\subsubsection{Derivatives of the Manufactured Solution}

Beginning with $\phi(x,t)$ the first and second derivatives of the manufactured solution with respect to time are:
\begin{align*}
  \phi_t(x,t) &= \frac{\partial}{\partial t} \Bigl[ 1 + \cos(\pi x)\cos(\pi t) \Bigr]
             = -\pi \cos(\pi x) \sin(\pi t), \\
  \phi_{tt}(x,t) &= \frac{\partial}{\partial t}\Bigl[-\pi \cos(\pi x) \sin(\pi t)\Bigr]
             = -\pi^2 \cos(\pi x) \cos(\pi t).
\end{align*}
Similarly, the spatial derivatives are:
\begin{align*}
  \phi_x(x,t) &= \frac{\partial}{\partial x} \Bigl[1 + \cos(\pi x) \cos(\pi t)\Bigr]
             = -\pi\sin(\pi x)\cos(\pi t), \\
  \phi_{xx}(x,t) &= -\pi^2 \cos(\pi x)\cos(\pi t).
\end{align*}

Thus, we observe that
\[
\phi_{tt}-\phi_{xx}= -\pi^2 \cos(\pi x)\cos(\pi t) - \left[-\pi^2 \cos(\pi x)\cos(\pi t)\right] = 0.
\]
Therefore, the left-hand side of the $\phi$-equation~\eqref{eq:phi_force} becomes
\[
\phi + V_\phi.
\]

Computing the time derivatives for $\chi(x,t)$:
\begin{align*}
  \chi_t(x,t) &= \frac{\partial}{\partial t} \Bigl[1 + \cos(\pi x) \sin(\pi t)\Bigr]
             = \pi \cos(\pi x) \cos(\pi t), \\
  \chi_{tt}(x,t) &= -\pi^2 \cos(\pi x) \sin(\pi t).
\end{align*}
And the spatial derivatives:
\begin{align*}
  \chi_x(x,t) &= \frac{\partial}{\partial x} \Bigl[1 + \cos(\pi x)\sin(\pi t)\Bigr]
             = -\pi \sin(\pi x) \sin(\pi t), \\
  \chi_{xx}(x,t) &= -\pi^2 \cos(\pi x) \sin(\pi t).
\end{align*}
Thus,
\[
\chi_{tt}-\chi_{xx} = -\pi^2 \cos(\pi x)\sin(\pi t) - \left[-\pi^2 \cos(\pi x)\sin(\pi t)\right] = 0.
\]
The left-hand side of the $\chi$-equation~\eqref{eq:chi_force} then becomes
\[
\chi - V_\chi.
\]

\subsubsection{Forcing Functions}
Thus, we define the forcing functions as:
\begin{align}
  F_\phi(x,t) &= \phi(x,t) + V_\phi(\phi(x,t),\chi(x,t)), \label{eq:Fphi}\\[1mm]
  F_\chi(x,t) &= \chi(x,t) - V_\chi(\phi(x,t),\chi(x,t)). \label{eq:Fchi}
\end{align}
With these definitions, the modified (forced) problem becomes:
\begin{align}
  \phi_{tt} - \phi_{xx} + \phi + V_\phi(\phi,\chi) &= F_\phi(x,t), \label{eq:phi_forced}\\[1mm]
  \chi_{tt} - \chi_{xx} + \chi - V_\chi(\phi,\chi) &= F_\chi(x,t). \label{eq:chi_forced}
\end{align}
By construction, substituting the manufactured solutions \eqref{eq:phi_manuf} and \eqref{eq:chi_manuf} into the left-hand sides of \eqref{eq:phi_forced} and \eqref{eq:chi_forced} yields exactly $F_\phi(x,t)$ and $F_\chi(x,t)$, respectively.

This recast's our weak EOM from Eq. \ref{eq:Weak_EOM}:
        \begin{align}
    K_1 &= \int\!\Bigl(u_t\Psi+\phi_x\Psi_x
             +m_\phi^2\phi\Psi \nonumber\\
        &\qquad+\partial_\phi V\,\Psi-F_\phi\Psi\Bigr),\nonumber\\
    K_2 &= \int(\phi_t\Psi-u\Psi),\nonumber\\
    G_1 &= \int\!\Bigl(v_t\Psi+\chi_x\Psi_x
             +m_\chi^2\chi\Psi \nonumber\\
        &\qquad-\partial_\chi V\,\Psi-F_\chi\Psi\Bigr),\nonumber\\
    G_2 &= \int(\chi_t\Psi-v\Psi).\nonumber
\end{align}
        
By incorporating the forcing functions into the numerical solver, we can compute numerical solutions and compare them to the exact manufactured solutions. A convergence study is shown in table \ref{tab:results_1+1Manufactured} which indicates 2nd order convergence.

\subsubsection{Results}

\begin{table}
  \centering
  \caption{Convergence Study Results (1+1)}
  \begin{tabular}{ccccc}
    \toprule
    Mesh ($nx\times nt$) & $\phi$ $L_{2}$ Error & Rate$_{\phi}$ & $\chi$ $L_{2}$ Error & Rate$_{\chi}$ \\ \midrule
    $100\times 101$   & 7.64e-04 & - & 9.08e-04 & -  \\[0.5em]
    $200\times 201$   & 1.90e-04 & 2.0037 & 2.26e-04 & 2.0049  \\[0.5em]
    $400\times 401$   & 4.76e-05 & 2.0018 & 5.65e-05 & 2.0023 \\[0.5em]
    $800\times 801$   & 1.19e-05 & 2.0008 & 1.41e-05 & 2.0011 \\ \bottomrule
  \end{tabular}
  \label{tab:results_1+1Manufactured}
\end{table}

The rates in table \ref{tab:results_1+1Manufactured} are calculated by doing 
\[
Rate_{\phi/\chi} = \frac{log\left(E_{coarse}/E_{fine}\right)}{log(h_{coarse}/h_{fine})}.
\]
The mesh is setup up so that $h = h_x = h_t$. Thus, $log(h_{coarse}/h_{fine}) = log(2)$.

Similarly, we tested a manufactured solution on the (2+1) code and obtained the results in table \ref{tab:results_2+1Manufactured}.

\begin{table}
  \centering
  \caption{Convergence Study Results (2+1)}
  \begin{tabular}{ccccc}
    \toprule
    Mesh ($nx\times ny\times nt$) & $\phi$ $L_{2}$ Error & Rate$_{\phi}$ & $\chi$ $L_{2}$ Error & Rate$_{\chi}$ \\ \midrule
    $10\times 10\times 11$   & 1.09e-01 & -    & 8.05e-02 & -    \\[0.5em]
    $20\times 20\times 21$   & 2.24e-02 & 2.286 & 1.73e-02 & 2.218 \\[0.5em]
    $40\times 40\times 41$   & 5.34e-03 & 2.071 & 4.21e-03 & 2.041 \\ \bottomrule
  \end{tabular}
  \label{tab:results_2+1Manufactured}
\end{table}

\begin{figure*}[!t]
    \centering
    \begin{subfigure}[t]{0.48\textwidth}
        \includegraphics[width=0.9\linewidth]{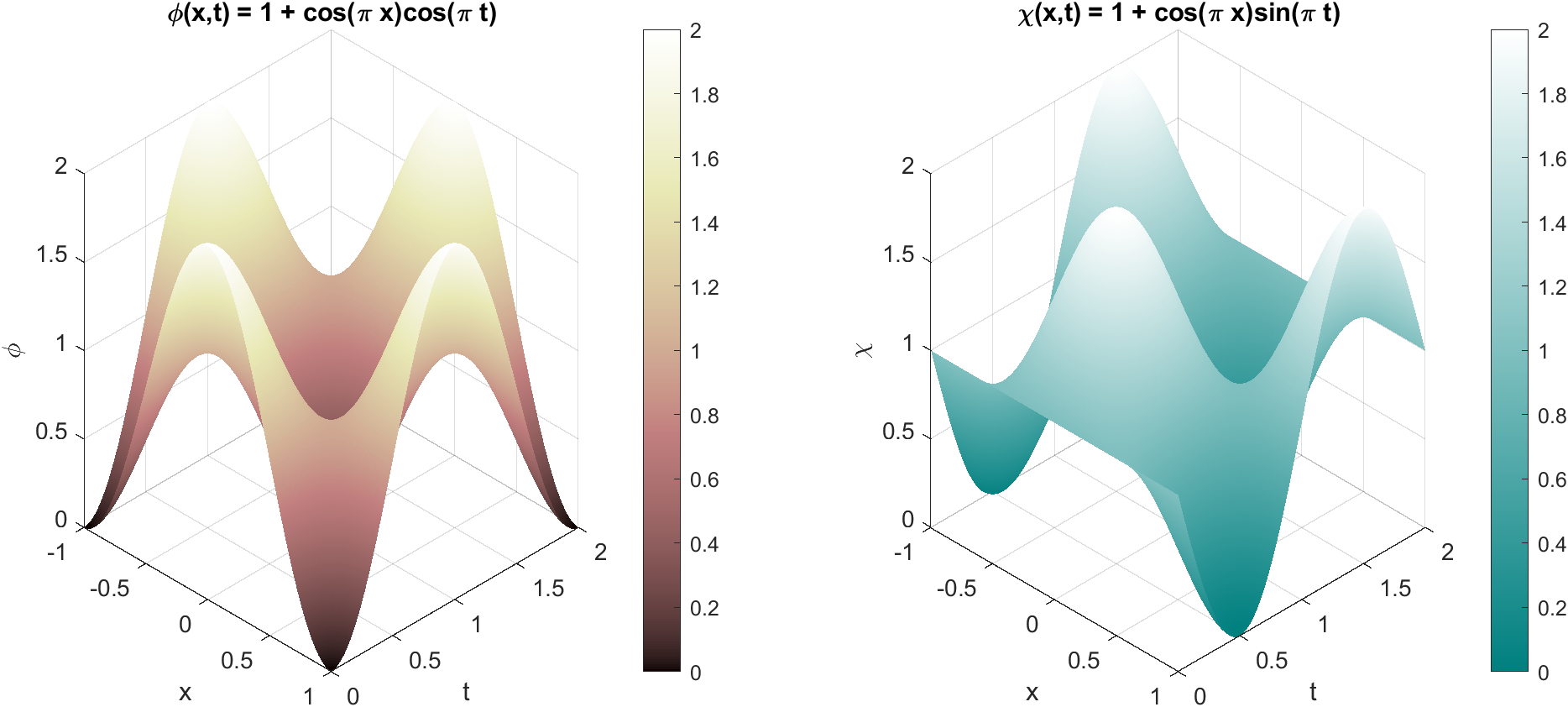}
        \caption{Analytic}
        \label{fig:MMS_analytic_3D}
    \end{subfigure}
    \hfill
    \begin{subfigure}[t]{0.48\textwidth}
        \includegraphics[width=0.9\linewidth]{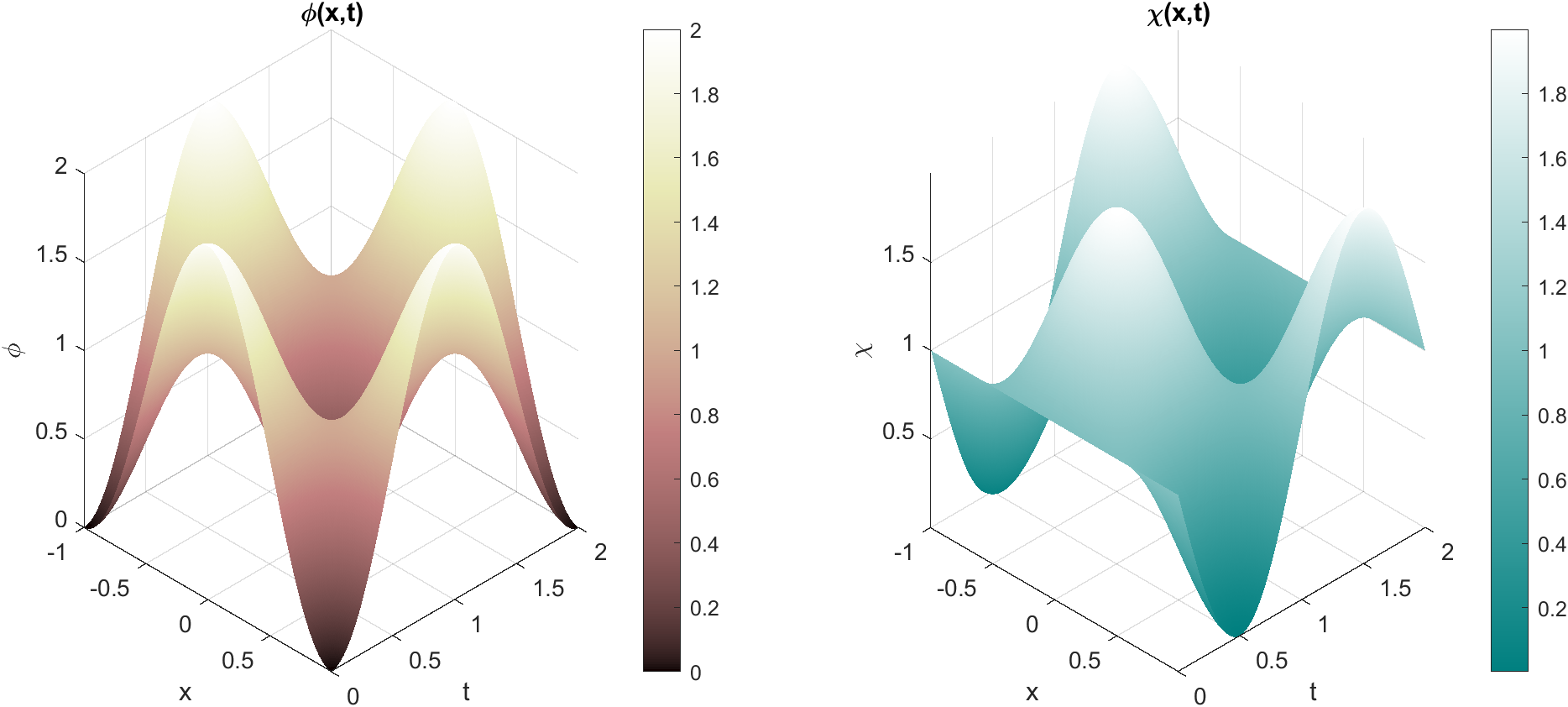}
        \caption{Numerical}
        \label{fig:MMS_numerical_3D}
    \end{subfigure}
    \caption{\justifying Comparison of spacetime surface plots for the analytic and numerical solutions in (1+1) dimensions.}
    \label{fig:MMS_3D_compare}
\end{figure*}

We see from figure \ref{fig:MMS_3D_compare} and table \ref{tab:results_1+1Manufactured} that the numerical results compare nicely with the manufactured results. These findings give us good reason to believe that our (1+1) numerical implementation is working properly.

\end{document}